\documentclass[11pt]{article}
\usepackage[utf8]{inputenc}
\usepackage[T1]{fontenc}
\usepackage{amsmath,amssymb,amsfonts,amsthm}
\usepackage{mathtools,lmodern}
\usepackage{geometry}
\usepackage[hidelinks,
pdftitle={Synchronized Normalized Solutions for a Critical Caffarelli-Kohn-Nirenberg System},
pdfauthor={Shengbing Deng and Yifan Zhao},
pdfkeywords={Caffarelli-Kohn-Nirenberg inequality; normalized solution; elliptic system; critical exponent; Pohozaev constraint}]{hyperref}
\newtheorem{theorem}{Theorem}[section]
\newtheorem{proposition}{Proposition}[section]
\newtheorem{lemma}{Lemma}[section]

\newtheorem{remark}{Remark}[section]

\numberwithin{equation}{section}

\title{Synchronized Normalized Solutions for a Critical Caffarelli--Kohn--Nirenberg System}
\author{Shengbing Deng\footnote{Corresponding author. E-mail address:\, {\tt shbdeng@swu.edu.cn}.}\ \ and Yifan Zhao\footnote{E-mail address:\, {\tt yfzhaomath@163.com}.}\\
	School of Mathematics and Statistics, Southwest University\\
	Chongqing, 400715, P.R. China}
\date{}

\begin{document}

\maketitle
\begin{abstract}
    We study positive synchronized normalized solutions to the critical
    Caffarelli--Kohn--Nirenberg system:
	\begin{equation*}
		\left\{
		\begin{aligned}
			-\operatorname{div}(|x|^{-2a}\nabla u)
			&=\lambda_1|x|^{-2a}u+\beta_1|x|^{-bq}|u|^{q-2}u\\
			&\quad+|x|^{-bp}\left(
			\mu_1|u|^{p-2}u
			+\nu\alpha|u|^{\alpha-2}u|v|^\gamma
			\right),\quad\text{in }\mathbb R^N,\\
			-\operatorname{div}(|x|^{-2a}\nabla v)
			&=\lambda_2|x|^{-2a}v+\beta_2|x|^{-bq}|v|^{q-2}v\\
			&\quad+|x|^{-bp}\left(
			\mu_2|v|^{p-2}v
			+\nu\gamma|u|^\alpha|v|^{\gamma-2}v
			\right),\quad\text{in }\mathbb R^N,
		\end{aligned}
		\right.
	\end{equation*}
	subject to the prescribed weighted masses
	\begin{equation*}
		\int_{\mathbb R^N}|x|^{-2a}u^2\,dx=\rho_1^2,
		\qquad
		\int_{\mathbb R^N}|x|^{-2a}v^2\,dx=\rho_2^2,
	\end{equation*}
	where $N\geq5$, $0\leq a<\frac{N-4}{2}$, $a<b<a+1$,
	$\alpha,\gamma>1$,
	$\alpha+\gamma=p=\frac{2N}{N-2(1+a-b)}$, $2<q<p$,
	$\beta_i,\mu_i,\rho_i>0$ for $i=1,2$, and $\nu>0$.
    The associated Lagrange multipliers
    $\lambda_1,\lambda_2\in\mathbb R$ are unknown.

    Under suitable algebraic compatibility conditions, we establish
    existence and multiplicity results for positive synchronized
    normalized solutions. The construction is based on synchronized
    reduction and critical energy estimates involving a finite-mass
    Caffarelli--Kohn--Nirenberg extremal.
    In the mass-subcritical regime, we find two positive solutions
    with opposite energies when
    $\beta=\beta_1\rho_1^{q-2}=\beta_2\rho_2^{q-2}$
    is sufficiently small. In the mass-critical and mass-supercritical
    regimes, we obtain a synchronized Pohozaev ground state below the
    critical energy threshold for sufficiently small $\beta$ and for all
    $\beta>0$, respectively. All these solutions satisfy $\lambda_1=\lambda_2<0$.
    For symmetric mass-critical data with unit weighted masses, we further
    show that the full vector Pohozaev level is strictly lower than the
    synchronized one under a weak coupling condition.
	

\end{abstract}

\section{Introduction and main results}\label{sec:introduction}

Coupled nonlinear Schrödinger equations arise in the study of interacting wave components, particularly in multicomponent Bose--Einstein condensates~\cite{Hall1998,Myatt1997} and nonlinear optics~\cite{ChristodoulidesJoseph1988}.
The interactions between different components play an
important role in the existence and structure of standing
waves \cite{BartschJeanjeanSoave2016,Sirakov2007}.
In normalized problems, the masses of the components
are prescribed, while the corresponding frequencies
are determined as Lagrange multipliers.
This formulation connects the search for standing-wave
profiles with critical points of the energy under mass
constraints.
For two components, the coupling must be treated
together with two prescribed masses and two unknown
multipliers.

The variational theory of normalized solutions for scalar nonlinear Schr\"odinger equations has attracted considerable attention. In the mass-subcritical regime, Cazenave
and Lions \cite{CazenaveLions1982} established the connection between
constrained minimization, existence, and orbital stability of standing
waves. When the energy is unbounded from below on the mass constraint,
direct minimization is no longer applicable. Jeanjean
\cite{Jeanjean1997} introduced a constrained mountain-pass framework
based on a scaling parameter and the Pohozaev identity, which yields
bounded Palais--Smale sequences.
Further developments in the variational theory of normalized
Schr\"odinger equations can be found in
\cite{AlvesJiMiyagaki2022,BieganowskiMederski2021,JeanjeanJendrejLeVisciglia2022,JeanjeanLu2020,Li2021,SoaveJDE2020}.

For problems involving combined nonlinearities with
Sobolev critical growth, a representative normalized model is
\[
\begin{cases}
	-\Delta u=\lambda u+\mu|u|^{q-2}u
	+|u|^{2^*-2}u, & \text{in }\mathbb R^N,\\[2mm]
	\displaystyle\int_{\mathbb R^N}|u|^2\,dx=\rho^2,
\end{cases}
\]
where $N\geq3$, $2<q<2^*=\frac{2N}{N-2}$,
$\mu,\rho>0$, and $\lambda\in\mathbb R$ is an unknown
Lagrange multiplier.
Soave \cite{SoaveJFA2020} investigated normalized
ground states for this type of problem and raised
the question of the existence of a second solution
of mountain-pass type in addition to the local minimizer.
In the mass-subcritical regime and under suitable
smallness assumptions, Jeanjean--Le \cite{JeanjeanLe2022}
and Wei--Wu \cite{WeiWu2022} subsequently established
the corresponding multiplicity results for
$N\geq4$ and $N=3$, respectively.
Related results for coupled Schr\"odinger systems,
including least-energy solutions and solutions
with prescribed masses, can be found in
\cite{BartschJeanjeanSoave2016,Deng2025,	Sirakov2007,YunZhang2022,Zhang2024}.
These developments naturally motivate the study of normalized
Schr\"odinger systems in weighted critical settings, where singular
weights introduce additional variational difficulties.

For weighted critical problems, the Caffarelli--Kohn--Nirenberg (CKN)
inequality provides a natural framework for the analysis of such
singular and critical nonlinearities.
The sharp CKN inequality takes the form
\begin{equation}\label{eq:sharp-ckn}
	S(a,b)\left(
	\int_{\mathbb R^N}|x|^{-bp}|w|^p\,dx
	\right)^{2/p}
	\leq
	\int_{\mathbb R^N}|x|^{-2a}|\nabla w|^2\,dx,
\end{equation}
where $S(a,b)$ denotes the best constant,
$p=\frac{2N}{N-2(1+a-b)}$,
$0\leq a<\frac{N-4}{2}$, and $a<b<a+1$.
The inequality was introduced by Caffarelli, Kohn
and Nirenberg \cite{CaffarelliKohnNirenberg1984}.
Catrina and Wang \cite{CatrinaWang2001}
investigated sharp constants and extremal functions,
while Felli and Schneider \cite{FelliSchneider2003}
established symmetry-breaking results.
Further developments concerning rigidity and stability
can be found in
\cite{DolbeaultEstebanLoss2016,FrankPeteranderl2024,
	WeiWuCKN2024,WeiWuCKN2026,ZhouZou2024}.
The condition $a\geq0$ places the sharp extremal
in the radial-symmetry range, while
$a<\frac{N-4}{2}$ ensures its square integrability
with respect to the weight $|x|^{-2a}$.
This finite-mass property is used in the critical energy estimates.

For normalized inhomogeneous Schrödinger equations with critical Hardy--Sobolev-type nonlinearities, Cardoso, de Oliveira and Miyagaki \cite{CardosoOliveiraMiyagaki2025} established related existence results. Their work provides a relevant scalar framework for comparison, particularly in the case \(a=0\), where the diffusion operator reduces to the standard Laplacian.
More recently, Goel and Rai \cite{GoelRai2026} studied a scalar
weighted equation with the CKN critical exponent in a radial setting.
Under suitable parameter restrictions, they obtained two positive
radial normalized solutions in the mass-subcritical regime and
established the existence of a positive radial normalized ground
state in the mass-critical and mass-supercritical regimes.
Their approach uses mass-preserving fiber maps to analyze the
Pohozaev geometry and construct constrained Palais--Smale sequences.
Strict energy estimates are then used to establish compactness
at the relevant variational levels.

The present study concerns the coupled CKN system, which we treat
through synchronized reduction.
The resulting scalar variational problem is analyzed in the full
weighted energy space of real-valued functions, without imposing
radial symmetry.
The main contributions of this paper concern the existence
and multiplicity of compatible synchronized solutions
to the coupled CKN system, together with a strict comparison
between the synchronized and full vector Pohozaev levels.

For prescribed masses $\rho_1,\rho_2>0$, we consider the system
\begin{equation}\label{eq:vector-system}
	\left\{
	\begin{aligned}
		-\operatorname{div}(|x|^{-2a}\nabla u)
		&=\lambda_1|x|^{-2a}u+\beta_1|x|^{-bq}|u|^{q-2}u\\
		&\quad+|x|^{-bp}\left(\mu_1|u|^{p-2}u
		+\nu\alpha |u|^{\alpha-2}u|v|^\gamma\right),\quad\text{in }\mathbb R^N,\\
		-\operatorname{div}(|x|^{-2a}\nabla v)
		&=\lambda_2|x|^{-2a}v+\beta_2|x|^{-bq}|v|^{q-2}v\\
		&\quad+|x|^{-bp}\left(\mu_2|v|^{p-2}v
		+\nu\gamma |u|^\alpha|v|^{\gamma-2}v\right),\quad\text{in }\mathbb R^N,\\
		\int_{\mathbb R^N}|x|^{-2a}u^2\,dx&=\rho_1^2,
		\qquad
		\int_{\mathbb R^N}|x|^{-2a}v^2\,dx=\rho_2^2,
	\end{aligned}
	\right.
\end{equation}
where $N\geq5$, $0\leq a<\frac{N-4}{2}$, $a<b<a+1$,
$d=1+a-b$, $p=\frac{2N}{N-2d}$, and $2<q<p$.
Here $\lambda_1,\lambda_2\in\mathbb R$ are unknown Lagrange multipliers.
The coupling exponents satisfy
\begin{equation}\label{eq:alpha-gamma-intro}
	\alpha>1,\qquad \gamma>1,\qquad \alpha+\gamma=p,
\end{equation}
and the coefficients satisfy
\begin{equation}\label{eq:positive-coefficients}
	\mu_1>0,\qquad \mu_2>0,\qquad \nu>0,\qquad
	\beta_1>0,\qquad \beta_2>0.
\end{equation}
To construct synchronized solutions, we impose the
compatibility conditions
\begin{equation}\label{eq:lower-compatibility}
	\beta_1\rho_1^{q-2}=\beta_2\rho_2^{q-2}=\beta>0
\end{equation}
and
\begin{equation}\label{eq:critical-compatibility}
	\mu_1\rho_1^{p-2}+\nu\alpha\rho_1^{\alpha-2}\rho_2^{\gamma}
	=\mu_2\rho_2^{p-2}+\nu\gamma\rho_1^{\alpha}\rho_2^{\gamma-2}
	=\kappa>0.
\end{equation}
Under these conditions, the synchronized ansatz
\[
u=\rho_1w,\qquad v=\rho_2w
\]
necessarily yields $\lambda_1=\lambda_2=:\lambda$.
Consequently, the system \eqref{eq:vector-system}
reduces, within the synchronized class, to the
scalar normalized problem
\begin{equation}\label{eq:scalar-equation-intro}
	\left\{
	\begin{aligned}
		-\operatorname{div}(|x|^{-2a}\nabla w)
		&=\lambda |x|^{-2a}w
		+\beta |x|^{-bq}|w|^{q-2}w
		+\kappa |x|^{-bp}|w|^{p-2}w,\\
		\int_{\mathbb R^N}|x|^{-2a}w^2\,dx&=1,
	\end{aligned}
	\right.
\end{equation}
where $\lambda\in\mathbb R$ is an unknown
Lagrange multiplier.

The main compactness difficulties arise from critical
concentration and the possible loss of prescribed masses
under weak convergence. Moreover, the vector Pohozaev
identity alone does not determine the signs of the two
Lagrange multipliers separately. The synchronized
reduction allows us to work with a single multiplier,
while recovering the prescribed mass and excluding
critical concentration remain essential.

Define the interpolation exponent and the mass-critical
exponent by
\[
\delta_q=\frac{(N-2a+2b)q-2N}{2q},
\qquad
q_c=\frac{2N+4}{N-2a+2b}.
\]
The exponent $q_c$ is determined by
the mass-preserving scaling \eqref{eq:mass-dilation-functional}.
Accordingly, $2<q<q_c$, $q=q_c$, and $q_c<q<p$ correspond to
the mass-subcritical, mass-critical, and mass-supercritical
regimes, respectively.

Let $E$ denote the energy functional associated with
\eqref{eq:vector-system}, defined in \eqref{eq:vector-energy}, and let $\mathcal S_{\rm syn}$
and $\mathcal P_{\rm syn}$ be the synchronized mass
constraint and Pohozaev set defined in
\eqref{eq:synchronized-sets}.
We denote by $K_q$ the optimal constant in the weighted
Gagliardo--Nirenberg inequality \eqref{eq:weighted-gn},
and by $T_\kappa$ the critical energy threshold for
the reduced scalar problem, defined in
Section~\ref{sec:functional}.
In the mass-subcritical regime,
$\mathcal P_{\rm syn}^{-}$ denotes the upper branch
of $\mathcal P_{\rm syn}$, while the positive smallness
threshold $\beta_0$ is specified in
Section~\ref{sec:fibers}. A synchronized Pohozaev ground state is understood as
a minimizer of $E$ over $\mathcal P_{\rm syn}$.
Our first main result is as follows.

\begin{theorem}\label{thm:main}
	Let $N\geq5,0\leq a<\frac{N-4}{2}, a<b<a+1,2<q<p$.
	Assume that \eqref{eq:alpha-gamma-intro}, \eqref{eq:positive-coefficients}, and the compatibility conditions \eqref{eq:lower-compatibility} and \eqref{eq:critical-compatibility} hold.
	
	\begin{enumerate}
		\item[(i)] If $2<q<q_c$, there exists
		\begin{equation*}
			\beta_0=\beta_0(N,a,b,q,\kappa)>0
		\end{equation*}
		such that, for $0<\beta<\beta_0$, system
		\eqref{eq:vector-system} admits two positive synchronized
		normalized solutions
		\begin{equation*}
			(u^+,v^+)=(\rho_1w^+,\rho_2w^+),
			\qquad
			(u^-,v^-)=(\rho_1w^-,\rho_2w^-).
		\end{equation*}
		The first is a local minimizer of $E$ on $\mathcal S_{\rm syn}$,
		while the second minimizes $E$ on $\mathcal P^-_{\rm syn}$.
		Their energies satisfy
		\begin{equation}\label{eq:opposite-energy-main}
			E(u^+,v^+)<0<E(u^-,v^-).
		\end{equation}
		
		\item[(ii)] Assume that either $q=q_c$ and
		\begin{equation}\label{eq:beta-critical-main}
			0<\beta<\frac{1}{\delta_qK_q},
		\end{equation}
		or $q_c<q<p$ and $\beta>0$ hold.
		Then system \eqref{eq:vector-system} admits a positive ground
		state on $\mathcal P_{\rm syn}$ of the form
		\begin{equation*}
			(u,v)=(\rho_1w,\rho_2w).
		\end{equation*}
		Its energy satisfies
		\begin{equation}\label{eq:critical-energy-main}
			0<E(u,v)<(\rho_1^2+\rho_2^2)T_{\kappa}.
		\end{equation}
	\end{enumerate}
	
	In all cases, the Lagrange multipliers satisfy
	\begin{equation}\label{eq:multiplier-main}
		\lambda_1=\lambda_2=\lambda
		=(\delta_q-1)\beta
		\int_{\mathbb R^N}|x|^{-bq}|w|^q\,dx<0.
	\end{equation}
\end{theorem}

\begin{remark}
	In the mass-critical case $q=q_c$, $\delta_q=2/q_c$. Thus, condition
	\eqref{eq:beta-critical-main} is equivalent to
	$\beta<q_c/(2K_q)$.
	For unit weighted mass and parameters in the common range,
	the corresponding threshold in Goel and Rai~\cite{GoelRai2026}
	has the same form, with their interpolation coefficient
	$C_{a,b}^{q_c}$ in place of $K_q$.
\end{remark}

\begin{remark}\label{rem:symmetric}
	The compatibility conditions \eqref{eq:lower-compatibility} and \eqref{eq:critical-compatibility} hold for the symmetric data $\rho_1=\rho_2=\rho$, $\beta_1=\beta_2$, $\mu_1=\mu_2$, and $\alpha=\gamma=p/2$, with
	$
    \beta=\beta_1\rho^{q-2},
    \kappa=\rho^{p-2}\left(\mu_1+\frac{p}{2}\nu\right).
    $
	
	More generally, these conditions also allow
	nonsymmetric data. With the masses, coupling exponents,
	and $\nu$ fixed, a family of compatible positive
	coefficients is given by
	\begin{equation*}
		\beta_i=\beta\rho_i^{2-q},\quad i=1,2,\qquad
		\mu_1=\rho_1^{2-p}
		\left(\kappa-\nu\alpha
		\rho_1^{\alpha-2}\rho_2^\gamma\right),\qquad
		\mu_2=\rho_2^{2-p}
		\left(\kappa-\nu\gamma
		\rho_1^\alpha\rho_2^{\gamma-2}\right),
	\end{equation*}
	where $\beta>0$ and
	$
	\kappa>\nu\max\left\{
	\alpha\rho_1^{\alpha-2}\rho_2^\gamma,\,
	\gamma\rho_1^\alpha\rho_2^{\gamma-2}
	\right\}.
	$
\end{remark}

The next result compares the synchronized and full vector Pohozaev
levels. The levels $m_{\rm vec}$, $m_{\rm syn}$, and $m_\chi$ are
defined in Section~\ref{sec:functional}.

\begin{theorem}[Strict vector comparison]\label{thm:weak-coupling}
Assume that \(N\geq5\), \(0\leq a<\frac{N-4}{2}\), \(a<b<a+1\), and \(q=q_c\), and that \eqref{eq:beta-critical-main} holds. Let
	\begin{equation*}
		\begin{aligned}
			&\rho_1=\rho_2=1,\qquad \beta_1=\beta_2=\beta,\qquad
			\mu_1=\mu_2=\mu>0,\\
			&\alpha=\gamma=p/2,\qquad \nu>0,\qquad
			\kappa=\mu+\frac p2\nu.
		\end{aligned}
	\end{equation*}
	Then
	\begin{equation}\label{eq:vector-level-comparison}
		0<m_{\rm vec}<m_\mu,
		\qquad m_{\rm syn}=2m_\kappa,
		\qquad m_\chi=\chi^{-2/(p-2)}m_1\quad(\chi>0).
	\end{equation}
	In particular, if
	\begin{equation}\label{eq:weak-coupling-interval}
		0<\nu\leq\frac{2\mu}{p}\left(2^{(p-2)/2}-1\right),
	\end{equation}
	then $m_{\rm vec}<m_{\rm syn}$. Every synchronized Pohozaev
	ground state provided by Theorem \ref{thm:main} has energy strictly greater than
	the full vector Pohozaev infimum.
\end{theorem}

The proof of Theorem~\ref{thm:main} combines synchronized
reduction with constrained variational methods.
The mass-preserving fiber maps characterize the Pohozaev
geometry and are used to construct bounded constrained
Palais--Smale sequences. Compactness at the relevant
energy levels is recovered through weighted subcritical
embeddings and an energy-splitting argument in the spirit
of Lions' concentration--compactness principle
\cite{Lions1984I,Lions1984II}. More precisely,
Brezis--Lieb splitting~\cite{BrezisLieb1983} and the sharp
CKN inequality yield the critical energy threshold
$T_\kappa$ for nontrivial concentration, while the
Lagrange multiplier identity recovers the prescribed
mass of any nontrivial weak limit. To exclude critical
concentration, we use a sharp CKN extremal with finite
weighted mass to derive strict energy estimates without
cutoff errors. In the mass-subcritical regime, the
lower-branch solution is obtained by local minimization,
while the upper-branch level lies strictly below the sum
of the lower-branch level and $T_\kappa$. In the
mass-critical and mass-supercritical regimes, the scalar
Pohozaev level lies strictly below $T_\kappa$.
These estimates yield the required compactness and
establish the scalar existence results, which are then
lifted to positive synchronized solutions of the
original system.

The proof of Theorem~\ref{thm:weak-coupling} is based on
a quantitative perturbation argument. We construct a
vector test pair whose first component is a scalar ground
state, while the second is a mass-normalized combination
of a small multiple of the same ground state and a diffuse
profile carrying the remaining mass. After mass normalization and a common
mass-preserving dilation, an energy expansion yields
$m_{\rm vec}<m_\mu$. The scaling law for $m_\chi$ and
the identity $m_{\rm syn}=2m_\kappa$ then imply
$m_{\rm vec}<m_{\rm syn}$ under the stated weak-coupling
condition. Attainment of $m_{\rm vec}$ requires a separate compactness
argument for both components and is not asserted here.

The remainder of the paper is organized as follows.
Section~\ref{sec:functional} introduces the functional
framework and establishes the synchronized reduction.
In Section~\ref{sec:fibers}, we analyze the geometry
of the fiber maps and construct constrained
Palais--Smale sequences.
Regularity properties, Pohozaev identities, and
compactness results are established in
Section~\ref{sec:compactness}.
The scalar existence results in the different mass
regimes are proved in Sections~\ref{sec:subcritical}
and~\ref{sec:critical-supercritical}.
Finally, Section~\ref{sec:main-proof} completes the
analysis of the vector system with a strict comparison
between the synchronized and full vector Pohozaev levels.

\section{Functional framework}\label{sec:functional}\label{sec:reduction}

We introduce the weighted functional framework and establish the
estimates and synchronized reduction used in the proofs.

For real-valued functions, we use the following notation
to facilitate subsequent computations:
\begin{align*}
	A(w)&=\int |x|^{-2a}|\nabla w|^2, & M(w)&=\int |x|^{-2a}w^2,\\
	B(w)&=\int |x|^{-bq}|w|^q, & C(w)&=\int |x|^{-bp}|w|^p.
\end{align*}
Let $D_a^{1,2}$ be the completion of $C_c^\infty(\mathbb R^N)$
with respect to the norm $A^{1/2}$, and set
\begin{equation*}
	X=D_a^{1,2}(\mathbb R^N)\cap L^2(\mathbb R^N;|x|^{-2a}),
	\quad \mathcal S=\{w\in X:M(w)=1\},
\end{equation*}
where $ \|w\|_X^2=A(w)+M(w).$
Unless otherwise specified, all spatial integrals
are taken over $\mathbb R^N$.

For $t\in\mathbb R$ and $w\in X$, define the
mass-preserving dilation
\begin{equation}\label{eq:mass-dilation-functional}
	(t\star w)(x)
	=e^{(N-2a)t/2}w(e^t x).
\end{equation}
A change of variables yields
\begin{equation}\label{eq:dilation-properties}
	\begin{aligned}
		M(t\star w)&=M(w),&
		A(t\star w)&=e^{2t}A(w),\\
		B(t\star w)&=e^{q\delta_q t}B(w),&
		C(t\star w)&=e^{pt}C(w).
	\end{aligned}
\end{equation}
In particular, $t\star w\in \mathcal S$ whenever $w\in \mathcal S$.

For comparison with the scalar formulation in
\cite{GoelRai2026}, the rescaling
$y=\kappa^{\frac{1}{p-2}}w,
\widehat\beta=\beta\kappa^{-\frac{q-2}{p-2}}
$
transforms \eqref{eq:scalar-equation-intro} into
an equivalent scalar problem with unit critical
coefficient, prescribed mass
$M(y)=\kappa^{2/(p-2)}$, and the same
Lagrange multiplier $\lambda$.
We retain $\kappa$ to make the synchronization
compatibility conditions explicit.

For later use, we record the following ground-state transform,
which relates the weighted space $X$ to the standard Sobolev space.

\begin{lemma}[Ground-state transform]\label{lem:ground-transform}
The map $w\mapsto |x|^{-a}w$ is an isomorphism from $X$ onto
$H^1(\mathbb R^N)$. For $z=|x|^{-a}w$, we have
\begin{equation}\label{eq:ground-state-identity}
 M(w)=\|z\|_2^2,
 \qquad A(w)=\int |\nabla z|^2-\mu_a\int |x|^{-2}z^2,
 \qquad \mu_a=a(N-2-a).
\end{equation}
Smooth functions compactly supported away from the origin are
dense in $X$.
\end{lemma}
\begin{proof}
	
For any $w\in C_c^\infty(\mathbb R^N\setminus\{0\})$,
using
\[
\operatorname{div}(|x|^{-2a-2}x)
=(N-2-2a)|x|^{-2a-2},
\]
integration by parts and the Hölder inequality yield
\begin{align*}
	(N-2-2a)\int_{\mathbb R^N}|x|^{-2a-2}w^2\,dx
	&=-2\int_{\mathbb R^N}
	|x|^{-2a-2}w\,x\cdot\nabla w\,dx\\
	&\leq 2\left(
	\int_{\mathbb R^N}|x|^{-2a-2}w^2\,dx
	\right)^{1/2}A(w)^{1/2}.
\end{align*}
Since $N-2-2a>0$, it follows that
\begin{equation}\label{eq:weighted-hardy}
    \int_{\mathbb R^N}|x|^{-2a-2}w^2\,dx
	\leq \frac{4}{(N-2-2a)^2}A(w).
\end{equation}
For any $\varphi\in C_c^\infty(\mathbb R^N)$,
let $\eta_\varepsilon$ be a smooth cutoff satisfying
$\eta_\varepsilon=0$ on $B_\varepsilon$,
$\eta_\varepsilon=1$ outside $B_{2\varepsilon}$,
and $|\nabla\eta_\varepsilon|\leq C/\varepsilon$.
Then
$\eta_\varepsilon\varphi\in
C_c^\infty(\mathbb R^N\setminus\{0\})$, and
\[
A(\eta_\varepsilon\varphi-\varphi)
\leq C_\varphi
\left(\varepsilon^{N-2a}
+\varepsilon^{N-2-2a}\right)
\longrightarrow0
\quad\text{as }\varepsilon\to0^+,
\]
since $N-2-2a>0$.
Consequently,
$C_c^\infty(\mathbb R^N\setminus\{0\})$
is dense in $D_a^{1,2}(\mathbb R^N)$,
and \eqref{eq:weighted-hardy} extends to
all $w\in D_a^{1,2}(\mathbb R^N)$
by continuity.

We next establish the density in $X$.
For any $w\in X$, choose a smooth radial cutoff
$\chi_{\varepsilon,R}$ such that
$\chi_{\varepsilon,R}=1$ on
$\{2\varepsilon\leq |x|\leq R\}$,
$\chi_{\varepsilon,R}=0$ on
$\{|x|\leq\varepsilon\}\cup\{|x|\geq2R\}$,
and
$|\nabla\chi_{\varepsilon,R}|\leq C/|x|$.
By the product rule,
\begin{align*}
	\|w-\chi_{\varepsilon,R}w\|_X^2
	&\leq C\int_{\{|x|<2\varepsilon\}\cup\{|x|>R\}}
	|x|^{-2a}(|\nabla w|^2+w^2)\,dx\\
	&\quad+C\int_{\{\varepsilon<|x|<2\varepsilon\}
		\cup\{R<|x|<2R\}}
	|x|^{-2a-2}w^2\,dx
	\longrightarrow0
\end{align*}
as $\varepsilon\to0^+$ and $R\to\infty$.
Indeed, the convergence follows from $w\in X$
and the integrability of
$|x|^{-2a-2}w^2$ ensured by
\eqref{eq:weighted-hardy}.
For fixed $\varepsilon$ and $R$,
the weights are bounded above and below
on the annulus containing the support of
$\chi_{\varepsilon,R}w$.
Standard mollification therefore yields
an approximation in
$C_c^\infty(\mathbb R^N\setminus\{0\})$
converging to $\chi_{\varepsilon,R}w$
in $X$.
Hence $C_c^\infty(\mathbb R^N\setminus\{0\})$
is dense in $X$.

For $w\in C_c^\infty(\mathbb R^N\setminus\{0\})$,
setting $z=|x|^{-a}w$, a direct computation
followed by integration by parts yields
\eqref{eq:ground-state-identity}.
By the classical Hardy inequality,
\[
\int_{\mathbb R^N}|x|^{-2}z^2\,dx
\leq \frac{4}{(N-2)^2}
\int_{\mathbb R^N}|\nabla z|^2\,dx,
\]
we obtain
\[
\left(1-\frac{4\mu_a}{(N-2)^2}\right)
\|\nabla z\|_2^2
\leq A(w)\leq\|\nabla z\|_2^2.
\]
Since
\[
1-\frac{4\mu_a}{(N-2)^2}
=\left(\frac{N-2-2a}{N-2}\right)^2>0,
\]
the gradient norms are equivalent.
Together with the mass identity
$M(w)=\|z\|_2^2$ and the density established
above, this shows that the transformation
$T:w\mapsto |x|^{-a}w$ extends to a
continuous linear map from $X$ into
$H^1(\mathbb R^N)$.

To establish surjectivity, let
$z\in H^1(\mathbb R^N)$.
Since $N\geq3$,
$C_c^\infty(\mathbb R^N\setminus\{0\})$
is dense in $H^1(\mathbb R^N)$.
Thus, there exists a sequence
$\{z_n\}\subset
C_c^\infty(\mathbb R^N\setminus\{0\})$
such that $z_n\to z$ in $H^1(\mathbb R^N)$.
Setting $w_n=|x|^a z_n$, the preceding
estimates imply that $\{w_n\}$ is Cauchy
in $X$. Hence $w_n\to w$ in $X$ for some
$w\in X$. Since
\[
\|w_n-|x|^a z\|_{L^2(\mathbb R^N;|x|^{-2a})}
=\|z_n-z\|_2\longrightarrow0,
\]
we have $w=|x|^a z$, and consequently
$Tw=z$. Moreover, the same estimates yield
\[
\|T^{-1}z\|_X^2
\leq \|z\|_{H^1(\mathbb R^N)}^2.
\]
Therefore, $T$ is a topological isomorphism
from $X$ onto $H^1(\mathbb R^N)$.
\end{proof}

The dilation \eqref{eq:mass-dilation-functional}
satisfies
\begin{equation}\label{eq:dilation-operator-bound}
	\|t\star h\|_X^2
	=e^{2t}A(h)+M(h)
	\leq e^{2|t|}\|h\|_X^2,
	\qquad t\in\mathbb R,\quad h\in X.
\end{equation}
By the density result in Lemma~\ref{lem:ground-transform}
and the estimate \eqref{eq:dilation-operator-bound},
the map $t\mapsto t\star h$ is continuous from
$\mathbb R$ into $X$ for every $h\in X$.

\begin{lemma}\label{lem:weighted-gn}
	For every $q\in(2,p]$, there exists an optimal constant
	$K_q>0$ such that
	\begin{equation}\label{eq:weighted-gn}
		\int_{\mathbb R^N}|x|^{-bq}|w|^q\,dx
		\leq K_q A(w)^{\frac{q\delta_q}{2}}
		M(w)^{\frac{q(1-\delta_q)}{2}},
		\qquad w\in X,
	\end{equation}
	where
	\[
	\delta_q=\frac{(N-2a+2b)q-2N}{2q}.
	\]
	In particular, $\delta_p=1$, and
	\eqref{eq:weighted-gn} coincides with the sharp
	Caffarelli--Kohn--Nirenberg inequality
	\eqref{eq:sharp-ckn} when $q=p$.
\end{lemma}

\begin{proof}
	For $q\in(2,p]$,
	\[
	\delta_q=b-a+N\left(\frac12-\frac1q\right),
	\qquad b-a<\delta_q\leq1.
	\]
	Moreover,
	\[
	\frac1q-\frac bN
	=\delta_q\left(\frac12-\frac{a+1}{N}\right)
	+(1-\delta_q)\left(\frac12-\frac aN\right).
	\]
	Setting
	$
	\sigma_q=-a-\frac{b-a}{\delta_q},
    $
	we obtain
	\[
	-b=\delta_q\sigma_q+(1-\delta_q)(-a),
	\qquad
	0<-a-\sigma_q=\frac{b-a}{\delta_q}<1.
	\]
	Since $a<(N-2)/2$ and $bq\leq bp<N$,
	the Caffarelli--Kohn--Nirenberg inequality
	\cite{CaffarelliKohnNirenberg1984} yields
	\begin{equation}\label{eq:ckn-interpolation-proof}
		\||x|^{-b}w\|_q
		\leq C_q
		\||x|^{-a}\nabla w\|_2^{\delta_q}
		\||x|^{-a}w\|_2^{1-\delta_q}
	\end{equation}
	for all $w\in C_c^\infty(\mathbb R^N)$.
	Raising both sides to the power $q$ and
	using the density of $C_c^\infty(\mathbb R^N)$
	in $X$, we obtain \eqref{eq:weighted-gn}.
	The optimal constant $K_q$ is therefore
	positive and finite.
\end{proof}

The functionals $B$ and $C$ are of class $C^1$ on $X$.
Indeed, for $r>2$ and real $s,t$,
\begin{equation*}
	\left||s|^{r-2}s-|t|^{r-2}t\right|
	\leq(r-1)(|s|+|t|)^{r-2}|s-t|.
\end{equation*}
Applying this scalar inequality to $|x|^{-b}w$ and using the H\"older inequality and Lemma~\ref{lem:weighted-gn}, we find that the derivatives of $B$ and $C$ are locally Lipschitz from $X$ to its dual. In particular, $A,B,C$ are uniformly Lipschitz on bounded subsets of $X$.

\begin{proposition}\label{prop:compact-embedding}
For every $q\in(2,p)$, the embedding
\begin{equation*}
 X\hookrightarrow L^q(\mathbb R^N;|x|^{-bq})
\end{equation*}
is compact.
\end{proposition}

\begin{proof}
Let $w_n\rightharpoonup w$ in $X$.  Put $z_n=w_n-w$.  Fix $0<R_0<R_1$.  The weights are bounded above and below on
\begin{equation*}
 \Omega_{R_0,R_1}=\{x:R_0<|x|<R_1\}.
\end{equation*}
Hence $z_n\rightharpoonup0$ in $H^1(\Omega_{R_0,R_1})$.  Since $p<2^*=\frac{2N}{N-2}$, Rellich's theorem gives
\begin{equation}\label{eq:annular-convergence}
 \int_{\Omega_{R_0,R_1}}|x|^{-bq}|z_n|^q\,dx\longrightarrow0.
\end{equation}
For $0<R<1$, the H\"older inequality and \eqref{eq:sharp-ckn} give
\begin{align}
 \int_{|x|<R}|x|^{-bq}|z_n|^q\,dx
 =\int_{|x|<R}\left(|x|^{-bp}|z_n|^p\right)^{\frac{q}{p}}\,dx\leq |B_R|^{1-\frac{q}{p}}C(z_n)^{\frac{q}{p}}
 \leq CR^{N\left(1-\frac{q}{p}\right)}.\label{eq:origin-tail}
\end{align}
Next, choose $\vartheta\in(0,1)$ such that
$
\frac1q=\frac{\vartheta}{2}
+\frac{1-\vartheta}{p}.
$
For $R>1$, the H\"older inequality yields
\begin{align}
 \int_{|x|>R}|x|^{-bq}|z_n|^q\,dx
 \leq R^{-\vartheta q(b-a)}
 M(z_n)^{\frac{\vartheta q}{2}}C(z_n)^{\frac{(1-\vartheta)q}{p}}\leq CR^{-\vartheta q(b-a)}.\label{eq:infinity-tail}
\end{align}
The bounds in \eqref{eq:origin-tail} and \eqref{eq:infinity-tail} are uniform in $n$.
Letting $R_0\to0^+$ and $R_1\to\infty$,
we conclude that
\[
w_n\longrightarrow w
\quad\text{in }L^q(\mathbb R^N;|x|^{-bq}).
\]
\end{proof}

The following lemma provides a sharp CKN extremal with finite
weighted mass. This property is essential for the subsequent
critical energy estimates.

\begin{lemma}\label{lem:finite-mass-extremal}
Let
$N\geq5, 0\leq a<\frac{N-4}{2}, a<b<a+1.
$
Then, the sharp constant $S(a,b)$ is attained by a positive radial function $W\in X$.  Up to multiplication and dilation, this function has the form
\begin{equation}\label{eq:explicit-extremal}
 W(x)=c_0\left(1+|x|^{\sigma}\right)^{-m},
 \qquad
 \sigma=\frac{2d(N-2-2a)}{N-2d},
 \qquad
 m=\frac{N-2d}{2d}.
\end{equation}
In particular,
\begin{equation*}
 \frac{A(W)}{C(W)^{2/p}}=S(a,b),
 \qquad
 M(W)<\infty.
\end{equation*}
\end{lemma}

\begin{proof}
The existence and radiality follow from the sharp Caffarelli--Kohn--Nirenberg theory in the range $a\geq0$; see \cite{CatrinaWang2001}.  The Emden--Fowler reduction gives \eqref{eq:explicit-extremal}.  The profile is bounded near the origin.  At infinity,
\begin{equation*}
 W(x)\asymp |x|^{-(N-2-2a)}.
\end{equation*}
Thus the weighted mass density is comparable to
\begin{equation*}
 r^{N-1-2a}r^{-2(N-2-2a)}=r^{-N+3+2a}.
\end{equation*}
It is integrable at infinity exactly when $a<(N-4)/2$. Hence $W\in X$.
\end{proof}

For the scalar problem \eqref{eq:scalar-equation-intro},
we define the critical energy threshold
\begin{equation*}
	T_{\kappa}
	=\frac{d}{N}S(a,b)^{\frac{N}{2d}}
	\kappa^{-\frac{N-2d}{2d}}
\end{equation*}
and the energy functional
\begin{equation*}
	J(w)=\frac12A(w)-\frac{\beta}{q}B(w)
	-\frac{\kappa}{p}C(w).
\end{equation*}
The corresponding Pohozaev functional is
\[
P(w)=\left.\frac{d}{dt}J(t\star w)\right|_{t=0}
=A(w)-\delta_q\beta B(w)-\kappa C(w),
\]
and the scalar Pohozaev set is
\[
\mathcal P=\{w\in\mathcal S:P(w)=0\}.
\]
The energy functional of the original coupled system is
\begin{align}
	E(u,v)
	&=\frac12\int_{\mathbb R^N}|x|^{-2a}
	\left(|\nabla u|^2+|\nabla v|^2\right)\,dx
	\notag-\frac1q\int_{\mathbb R^N}|x|^{-bq}
	\left(\beta_1|u|^q+\beta_2|v|^q\right)\,dx
	\notag\\
	&\quad-\frac1p\int_{\mathbb R^N}|x|^{-bp}
	\left(\mu_1|u|^p+\mu_2|v|^p
	+p\nu|u|^\alpha|v|^\gamma\right)\,dx.
	\label{eq:vector-energy}
\end{align}
We denote the mixed nonlinear term by
\[
D(u,v)=\int_{\mathbb R^N}
|x|^{-bp}|u|^\alpha|v|^\gamma\,dx.
\]
By the H\"older inequality,
$
D(u,v)\leq C(u)^{\alpha/p}C(v)^{\gamma/p}.
$
Since $\alpha,\gamma>1$, the weighted
Sobolev embeddings imply that
$E\in C^1(X^2,\mathbb R)$.
The corresponding mass constraint and
Pohozaev set are defined by
\[
\mathcal S_{\rm vec}
=\{(u,v)\in X^2:
M(u)=\rho_1^2,\ M(v)=\rho_2^2\},
\]
and
\[
\mathcal P_{\rm vec}
=\left\{(u,v)\in\mathcal S_{\rm vec}:
\left.\frac{d}{dt}E(t\star u,t\star v)
\right|_{t=0}=0\right\},
\]
respectively.

Finally, we introduce the synchronized
mass constraint and Pohozaev set
\begin{equation}\label{eq:synchronized-sets}
		\mathcal S_{\rm syn}
		=\{(\rho_1w,\rho_2w):w\in\mathcal S\},\quad
		\mathcal P_{\rm syn}
		=\{(\rho_1w,\rho_2w):w\in\mathcal P\}.
\end{equation}
In the mass-subcritical regime, we further define
\[
\mathcal P_{\rm syn}^-
=\{(\rho_1w,\rho_2w):w\in\mathcal P^-\},
\]
where $\mathcal P^-$ denotes the upper scalar
Pohozaev branch characterized in
Lemma \ref{lem:subcritical-fibers}.
Set
\begin{equation*}
	m_{\rm vec}=\inf_{\mathcal P_{\rm vec}}E,
	\qquad m_{\rm syn}=\inf_{\mathcal P_{\rm syn}}E.
\end{equation*}
For $\chi>0$, write
\begin{equation*}
	J_\chi=\frac12A-\frac\beta qB-\frac\chi pC,
	\qquad
	\mathcal P_\chi=\{w\in\mathcal S:A(w)-\delta_q\beta B(w)-\chi C(w)=0\},
\end{equation*}
and let $m_\chi=\inf_{\mathcal P_\chi}J_\chi$.

\begin{proposition}\label{prop:exact-reduction}
	Assume that \eqref{eq:lower-compatibility} and
	\eqref{eq:critical-compatibility} hold.
	If $w\in\mathcal S$ is a constrained critical
	point of $J$, then
	\[
	(u,v)=(\rho_1w,\rho_2w)
	\]
	is a normalized solution of
	\eqref{eq:vector-system} with
	$\lambda_1=\lambda_2=\lambda$.
	Moreover,
	\begin{equation}\label{eq:energy-factorization}
		E(\rho_1w,\rho_2w)
		=(\rho_1^2+\rho_2^2)J(w).
	\end{equation}
\end{proposition}

\begin{proof}	
	Since $w$ is a constrained critical point
	of $J$ on $\mathcal S$, it satisfies
	\eqref{eq:scalar-equation-intro} for some
	$\lambda\in\mathbb R$.
	The compatibility conditions \eqref{eq:lower-compatibility}
	and \eqref{eq:critical-compatibility} give
	\[
	\beta\rho_i=\beta_i\rho_i^{q-1},
	\qquad i=1,2,
	\]
	and
	\[
	\kappa\rho_1
	=\mu_1\rho_1^{p-1}
	+\nu\alpha\rho_1^{\alpha-1}\rho_2^\gamma,\qquad
	\kappa\rho_2
	=\mu_2\rho_2^{p-1}
	+\nu\gamma\rho_1^\alpha\rho_2^{\gamma-1}.
	\]
	Multiplying \eqref{eq:scalar-equation-intro}
	by $\rho_1$ and $\rho_2$, respectively,
	and using these identities, we obtain
	\eqref{eq:vector-system} with
	$\lambda_1=\lambda_2=\lambda$.
	The prescribed mass constraints follow
	from $M(w)=1$.
	
	To prove \eqref{eq:energy-factorization},
	set $R=\rho_1^2+\rho_2^2$.
	The compatibility conditions yield
	\[
	\beta_1\rho_1^q+\beta_2\rho_2^q
	=\beta R
	\]
	and, since $\alpha+\gamma=p$,
	\[
	\mu_1\rho_1^p+\mu_2\rho_2^p
	+p\nu\rho_1^\alpha\rho_2^\gamma
	=\kappa R.
	\]
	Consequently, by \eqref{eq:vector-energy},
	\[
		E(\rho_1w,\rho_2w)
		=\frac{R}{2}A(w)
		-\frac{\beta R}{q}B(w)
		-\frac{\kappa R}{p}C(w)
		=RJ(w).
	\]
\end{proof}

The Pohozaev functional associated with
the common mass-preserving dilation is
given by
\[
\begin{aligned}
	P_{\rm vec}(u,v)
	&=A(u)+A(v)
	-\delta_q\bigl(\beta_1B(u)+\beta_2B(v)\bigr)\\
	&\quad-\mu_1C(u)-\mu_2C(v)
	-p\nu D(u,v).
\end{aligned}
\]
Indeed,
\[
P_{\rm vec}(u,v)
=\left.\frac{d}{dt}
E(t\star u,t\star v)\right|_{t=0}.
\]
By the compatibility conditions and
\eqref{eq:energy-factorization}, we obtain
\[
P_{\rm vec}(\rho_1w,\rho_2w)
=(\rho_1^2+\rho_2^2)P(w),
\qquad w\in\mathcal S.
\]
Consequently,
\[
\mathcal P_{\rm syn}
=\mathcal P_{\rm vec}\cap\mathcal S_{\rm syn}.
\]
Moreover, \eqref{eq:energy-factorization}
shows that minimization of $J$ over
$\mathcal P$ is equivalent to minimization
of $E$ over $\mathcal P_{\rm syn}$.
This correspondence is restricted to the synchronized class
and does not characterize general vector solutions.

\section{Fiber maps and constrained Palais--Smale sequences}\label{sec:fibers}

In this section, we use the framework of
Section~\ref{sec:functional} to establish the fiber geometry
in Lemmas~\ref{lem:subcritical-fibers}
and~\ref{lem:critical-supercritical-fibers}.
We then develop the tools for constructing constrained
Palais--Smale sequences.

For $w\in \mathcal S$, the mass-preserving dilation
\eqref{eq:mass-dilation-functional} and
\eqref{eq:dilation-properties} give the fiber map
\[
\Phi_w(t)=J(t\star w)
=\frac12 e^{2t}A(w)
-\frac{\beta}{q}e^{q\delta_q t}B(w)
-\frac{\kappa}{p}e^{pt}C(w).
\]
Writing $\theta=q\delta_q$, its derivatives are
\begin{equation}\label{eq:fiber-first}
 \Phi_w'(t)=e^{2t}A(w)-\delta_q\beta e^{\theta t}B(w)
 -\kappa e^{pt}C(w)
\end{equation}
and
\begin{equation*}
 \Phi_w''(t)=2e^{2t}A(w)-q\delta_q^2\beta e^{\theta t}B(w)
 -p\kappa e^{pt}C(w).
\end{equation*}
If $\theta<2$, we partition the Pohozaev set $\mathcal P$ according
to the sign of $\Phi_w''(0)$ as
\begin{equation*}
	\mathcal P=\mathcal P^+\cup\mathcal P^-\cup\mathcal P^0,
\end{equation*}
where
\begin{align*}
	\mathcal P^+
	&=\bigl\{w\in\mathcal P:\Phi_w''(0)>0\bigr\},\\
	\mathcal P^-
	&=\bigl\{w\in\mathcal P:\Phi_w''(0)<0\bigr\},\\
	\mathcal P^0
	&=\bigl\{w\in\mathcal P:\Phi_w''(0)=0\bigr\}.
\end{align*}

On $\mathcal P$,
\begin{equation}\label{eq:energy-on-P-A}
 J(w)=\frac{d}{N}A(w)
 -\left(\frac1q-\frac{\delta_q}{p}\right)\beta B(w)
\end{equation}
and
\begin{equation}\label{eq:energy-on-P-C}
 J(w)=\left(\frac{\delta_q}{2}-\frac1q\right)\beta B(w)
 +\frac{d}{N}\kappa C(w).
\end{equation}
For $s\geq0$, define
\begin{equation*}
 f_{\beta}(s)=\frac12s^2-\frac{\beta K_q}{q}s^{\theta}
 -\frac{\kappa}{pS(a,b)^{p/2}}s^p.
\end{equation*}

In the mass-subcritical case,
$0<\theta<2$. To specify the smallness constant used below, set
\begin{equation*}
 c_*=\frac{p-2}{p-\theta}
 \left(\frac{(2-\theta)S(a,b)^{\frac{p}{2}}}{(p-\theta)\kappa}\right)^{\frac{2-\theta}{p-2}}
\end{equation*}
and
\begin{equation*}
 \beta_{\rm f}=\frac{c_*}{2\delta_qK_q}.
\end{equation*}
Choose
\begin{equation}\label{eq:R-star}
 R_*=\frac12
 \left(\frac{pS(a,b)^{p/2}}{8\kappa}\right)^{1/(p-2)}
\end{equation}
and set
\begin{equation}\label{eq:beta-barrier}
 \beta_{\rm b}=\frac{q}{8K_q}R_*^{2-\theta}.
\end{equation}
We use the constant
\begin{equation*}
 \beta_0=\min\{\beta_{\rm f},\beta_{\rm b}\}.
\end{equation*}
For the mass-critical case, assumption \eqref{eq:beta-critical-main} ensures that
\[
c_\beta:=1-\delta_q\beta K_q>0.
\]

\begin{lemma}\label{lem:subcritical-fibers}
Let $2<q<q_c$ and $0<\beta< \beta_0$. Then $\mathcal P^0=\varnothing$.  Each fiber map has exactly two critical points
\begin{equation*}
 t_1(w)<t_2(w).
\end{equation*}
The first is a strict local minimum with negative energy, and the second is a strict global maximum with positive energy.  Moreover,
\begin{equation*}
 t_1(w)\star w\in\mathcal P^+,
 \qquad
 t_2(w)\star w\in\mathcal P^-.
\end{equation*}
Both maps $w\mapsto t_i(w)$ are of class $C^1$.
\end{lemma}

\begin{proof}
Put $s=e^t$.  Equation \eqref{eq:fiber-first} is equivalent to
\begin{equation}\label{eq:g-equation}
 g(s)=\delta_q\beta B(w),
 \qquad
 g(s)=A(w)s^{2-\theta}-\kappa C(w)s^{p-\theta}.
\end{equation}
Direct differentiation shows that $g$ has a unique maximum, given by
\begin{align*}
	\max_{s>0}g(s)
	&=\frac{p-2}{p-\theta}A(w)
	\left(\frac{(2-\theta)A(w)}{(p-\theta)\kappa C(w)}\right)^{\frac{2-\theta}{p-2}}.
\end{align*}
By \eqref{eq:sharp-ckn},
\begin{equation*}
 C(w)\leq S(a,b)^{-p/2}A(w)^{p/2}.
\end{equation*}
Therefore
\begin{equation*}
 \max_{s>0}g(s)\geq c_*A(w)^{\theta/2}.
\end{equation*}
Since $M(w)=1$, \eqref{eq:weighted-gn} gives
\begin{equation}\label{eq:B-upper-fiber}
 \delta_q\beta B(w)\leq\delta_q\beta K_qA(w)^{\theta/2}
 <\frac12c_*A(w)^{\theta/2}.
\end{equation}
Thus \eqref{eq:g-equation} has exactly two positive roots.  The sign pattern of $\Phi_w'$ gives the local minimum and the global maximum.  Since $\theta<2$,
\begin{equation*}
 \Phi_w(t)=-\frac{\beta}{q}e^{\theta t}B(w)+o(e^{\theta t})<0
 \quad\text{as }t\to-\infty.
\end{equation*}
Hence the first critical value is negative.

For the second value, applying \eqref{eq:weighted-gn} and \eqref{eq:sharp-ckn} yields
\begin{equation*}
 \Phi_w(t)\geq\frac12y^2-\frac{\beta K_q}{q}y^{\theta}
 -\frac{\kappa}{pS(a,b)^{p/2}}y^p,
 \qquad y=e^tA(w)^{1/2}.
\end{equation*}
Taking $y=R_*$ and using \eqref{eq:beta-barrier}, we obtain
\begin{equation}\label{eq:fiber-positive-at-R}
 \Phi_w(t_R)\geq f_{\beta}(R_*)\geq\frac14R_*^2>0,
 \qquad e^{t_R}A(w)^{1/2}=R_*.
\end{equation}
Thus the global maximum is positive.

The strict inequality in \eqref{eq:B-upper-fiber}
rules out double roots of \eqref{eq:g-equation},
and hence $\mathcal P^0=\varnothing$.
The final assertion follows from the implicit
function theorem.
\end{proof}

\begin{lemma}\label{lem:critical-supercritical-fibers}
Let $q_c\leq q<p$. If $q=q_c$, assume in addition that
\eqref{eq:beta-critical-main} holds. Then, for every
$w\in\mathcal S$, the fiber map $\Phi_w$ has exactly one
critical point $t_*(w)$, which is a strict global maximum.
Moreover, the map $w\mapsto t_*(w)$ is of class $C^1$ on
$\mathcal S$.
\end{lemma}

\begin{proof}
If $q>q_c$, then $\theta>2$.  Dividing \eqref{eq:fiber-first} by $e^{2t}$ gives
\begin{equation*}
 A(w)=\delta_q\beta e^{(\theta-2)t}B(w)
 +\kappa e^{(p-2)t}C(w).
\end{equation*}
The right-hand side increases strictly from $0$ to $+\infty$, so the equation has a unique root.

If $q=q_c$, then $\theta=2$.  By \eqref{eq:weighted-gn},
\begin{equation*}
 A(w)-\delta_q\beta B(w)
 \geq c_\beta A(w)>0.
\end{equation*}
Equation \eqref{eq:fiber-first} therefore admits a unique
root $t_*(w)$. In both cases, $\Phi_w'(t)$ changes sign
from positive to negative at $t_*(w)$, so that
$t_*(w)$ is the unique strict global maximizer of $\Phi_w$.
Moreover, the identity $\Phi_w'(t_*(w))=0$ gives
\[
	\Phi_w''(t_*(w))
	={}-\delta_q\beta(\theta-2)e^{\theta t_*(w)}B(w)-(p-2)\kappa e^{pt_*(w)}C(w)<0.
\]
The implicit function theorem then implies that
the map $w\mapsto t_*(w)$ is of class $C^1$ on $\mathcal S$.
\end{proof}

We now establish the variational properties common to
the three mass regimes. Throughout the remainder of this
section, we set
\[
\mathcal H=
\begin{cases}
	\mathcal P^-, & 2<q<q_c,\\
	\mathcal P, & q_c\le q<p.
\end{cases}
\]
We denote by $t_*(w)$ the unique global maximizer of
$\Phi_w$, with $t_*(w)=t_2(w)$ when $q<q_c$.
We assume $\beta<\beta_0$ in the mass-subcritical case
and \eqref{eq:beta-critical-main} in the mass-critical
case.

\begin{lemma}[Uniform fiber gap]\label{lem:fiber-gap}
	There exists a constant $c_1>0$ such that
	\begin{equation}\label{eq:uniform-fiber-gap}
	 -\Phi_z''(0)\geq c_1\qquad(z\in\mathcal H).
	\end{equation}
	Moreover, for every $L>0$, there exist constants
	$r_L,C_L>0$ such that
	\begin{equation}\label{eq:uniform-fiber-projection}
		|t_*(y)|\leq C_L\|y-z\|_X
	\end{equation}
	whenever $z\in\mathcal H$, $\|z\|_X\leq L$,
	$y\in\mathcal S$, and $\|y-z\|_X<r_L$.
\end{lemma}

\begin{proof}
    Suppose $\theta<2$. For $z\in\mathcal P^-$, the conditions
    $P(z)=0$ and $\Phi_z''(0)<0$ imply that
	\[
	A(z)<\frac{p-\theta}{2-\theta}\kappa C(z).
	\]
	Together with \eqref{eq:sharp-ckn}, this implies
	\[
	A(z)\geq A_*:=
	\left(
	\frac{(2-\theta)S(a,b)^{p/2}}
	{(p-\theta)\kappa}
	\right)^{2/(p-2)}.
	\]
	Since
	$
	c_*=\frac{p-2}{p-\theta}A_*^{(2-\theta)/2},
	$
	the assumption $\beta<\beta_{\rm f}$, together
	with Lemma~\ref{lem:weighted-gn}, gives
	\[
	\begin{aligned}
		-\Phi_z''(0)
		&= (p-2)A(z)
		-\delta_q(p-\theta)\beta B(z)\\
		&\geq
		\left(
		p-2-\delta_q(p-\theta)\beta K_q
		A_*^{(\theta-2)/2}
		\right)A(z)\\
		&\geq \frac{p-2}{2}A(z)
		\geq \frac{p-2}{2}A_*>0.
	\end{aligned}
	\]
    Next, suppose $\theta>2$. The identity $P(z)=0$ yields
	\[
		-\Phi_z''(0)
		=\delta_q(\theta-2)\beta B(z)
		+(p-2)\kappa C(z)\geq
		\min\{\theta-2,p-2\}A(z).
	\]
	By Lemma~\ref{lem:weighted-gn} and
	\eqref{eq:sharp-ckn}, we obtain
	\[
	1\leq
	\delta_q\beta K_q A(z)^{(\theta-2)/2}
	+\kappa S(a,b)^{-p/2}A(z)^{(p-2)/2}.
	\]
	Since $\theta>2$ and $p>2$, this inequality
	implies that $A(z)\geq A_0>0$ for some
	constant $A_0$ independent of $z$. Consequently, $-\Phi_z''(0)$ admits a
	uniform positive lower bound.
	Finally, suppose $\theta=2$. Recalling that
	$c_\beta=1-\delta_q\beta K_q>0$, we have
	\[
	-\Phi_z''(0)
	=(p-2)\kappa C(z)
	\geq(p-2)c_\beta A(z).
	\]
	Moreover, $P(z)=0$ and
	\eqref{eq:sharp-ckn} imply
	$
	A(z)\geq
	\left(
	\frac{c_\beta S(a,b)^{p/2}}{\kappa}
	\right)^{\frac{2}{p-2}}.
	$
	Thus \eqref{eq:uniform-fiber-gap} holds
	in all three cases.
	
	We next prove \eqref{eq:uniform-fiber-projection}.
	Define
	\[
	H(t,y)=\Phi_y'(t)
	=e^{2t}A(y)
	-\delta_q\beta e^{\theta t}B(y)
	-\kappa e^{pt}C(y).
	\]
	Fix $L>0$ and let $z\in\mathcal H$ satisfy
	$\|z\|_X\leq L$. Since $t_*(z)=0$, we have
	\[
	H(0,z)=0,
	\qquad
	H_t(0,z)=\Phi_z''(0)\leq-c_1.
	\]
	Choose $r_L\leq1$ so that $\|y\|_X\leq L+1$
	whenever $\|y-z\|_X<r_L$.
Since $A$, $B$, and $C$ are uniformly Lipschitz
on bounded subsets of $X$, and $H(0,z)=0$,
there exists a constant $K_L>0$, independent
of $y$ and $z$, such that
\[
|H(0,y)|+|H_t(0,y)-H_t(0,z)|
\leq K_L\|y-z\|_X.
\]
Furthermore, the explicit expression
for $H$ yields a constant $M_L>0$ such that
\[
|H_{tt}(t,y)|\leq M_L
\]
for $|t|\leq1$ and $\|y\|_X\leq L+1$.
Consequently,
\[
H_t(t,y)
\leq -c_1+M_L|t|+K_L\|y-z\|_X.
\]
	Choose $h_L,r_L\in(0,1)$ sufficiently
	small so that
	$
	M_Lh_L+K_Lr_L\leq\frac{c_1}{2},
	K_Lr_L<\frac{c_1h_L}{4}.
	$
	Then, whenever $\|y-z\|_X<r_L$,
	\[
	H_t(t,y)\leq-\frac{c_1}{2}
	\qquad (|t|\leq h_L).
	\]
	It follows that
	\[
	H(-h_L,y)>0,
	\qquad
	H(h_L,y)<0.
	\]
	Therefore, $H(\cdot,y)$ has a unique zero
	in $(-h_L,h_L)$. By the fiber geometry
	established in Lemma~\ref{lem:subcritical-fibers} and Lemma~\ref{lem:critical-supercritical-fibers},
	this zero coincides with $t_*(y)$.
	Finally, the mean value theorem yields
	\[
	\frac{c_1}{2}|t_*(y)|
	\leq |H(0,y)|
	\leq K_L\|y-z\|_X.
	\]
	Thus \eqref{eq:uniform-fiber-projection}
	follows with $C_L=2K_L/c_1$.
\end{proof}

\begin{lemma}[Maximum functional]\label{lem:maximum-functional}
The functional
\begin{equation}\label{eq:maximum-functional}
 F(w)=\max_{t\in\mathbb R}J(t\star w)
 =\frac12e^{2t_*(w)}A(w)-\frac\beta q e^{\theta t_*(w)}B(w)
 -\frac\kappa p e^{pt_*(w)}C(w)
\end{equation}
is of class $C^1$ on $\mathcal S$, and its derivative is given by
\begin{equation}\label{eq:maximum-derivative}
 dF(w)[h]=J'(t_*(w)\star w)[t_*(w)\star h]
 \quad(h\in T_w\mathcal S).
\end{equation}
Every local minimizer of $J$ on $\mathcal H$ is a constrained
critical point of $J$ on $\mathcal S$.
\end{lemma}

\begin{proof}
Differentiate \eqref{eq:maximum-functional}. The coefficient
of $dt_*(w)[h]$ vanishes because $H(t_*(w),w)=0$. The remaining terms yield \eqref{eq:maximum-derivative}, since the dilation is bounded and linear on $X$.
Since $t_*$ is continuous on $\mathcal S$,
the strong continuity of the dilation
and \eqref{eq:dilation-operator-bound}
imply that the map
$w\mapsto t_*(w)\star w$
is continuous from $\mathcal S$ into $\mathcal H$.
Let $z\in\mathcal H$ be a local minimizer of
$J$ on $\mathcal H$. Since $t_*(z)=0$,
it follows that $z$ is a local minimizer of
$F$ on $\mathcal S$. Consequently,
$dF(z)=0$, and \eqref{eq:maximum-derivative}
yields
$
d(J|_{\mathcal S})(z)=0.
$
\end{proof}

\begin{proposition}[Constrained Palais--Smale sequence]\label{prop:full-PS}
Let $c=\inf_{\mathcal H}J$. If there exists a bounded minimizing
sequence for $c$ in $\mathcal H$, then there exists a bounded Palais--Smale sequence
$\{w_n\}\subset\mathcal H$ for $J|_{\mathcal S}$
at level $c$ satisfying
\begin{equation}\label{eq:full-PS-sequence}
	J(w_n)\to c,\qquad P(w_n)=0,
	\qquad \|d(J|_{\mathcal S})(w_n)\|\to0.
\end{equation}
\end{proposition}

\begin{proof}
The fiber lemmas give $F>0$ and
$\inf_{\mathcal S}F=\inf_{\mathcal H}J=c$.
	Let $\{z_n\}\subset\mathcal H$ be a bounded
	minimizing sequence such that
	\[
	J(z_n)\leq c+\varepsilon_n^2,
	\qquad \varepsilon_n\to0^+.
	\]
Since $\mathcal S$ is closed in $X$, it is complete with respect to the induced metric. Applying Ekeland's variational
	principle to $F$ on $\mathcal S$, we obtain
	$\{y_n\}\subset\mathcal S$ satisfying
	\[
	\|y_n-z_n\|_X\leq\varepsilon_n,
	\quad
	c\leq F(y_n)\leq F(z_n),\quad\|dF(y_n)\|_{(T_{y_n}\mathcal S)^*}
	\leq\varepsilon_n.
	\]
The derivative estimate follows by applying
Ekeland's variational inequality along the
normalized curves
$
\gamma(t)=
\frac{y_n+th}{M(y_n+th)^{1/2}},
h\in T_{y_n}\mathcal S,
$
which satisfy $\gamma(0)=y_n$ and
$\gamma'(0)=h$.
	
	Since $\{z_n\}$ is bounded, there exists
	$L>0$ such that $\|z_n\|_X\leq L$ for all $n$.
	By \eqref{eq:uniform-fiber-projection},
	for all sufficiently large $n$,
	\[
	|s_n|:=|t_*(y_n)|
	\leq C_L\|y_n-z_n\|_X
	\leq C_L\varepsilon_n\to0.
	\]
	Set $w_n=s_n\star y_n$. Then
	$w_n\in\mathcal H$ and
	\[
	J(w_n)=F(y_n)\to c,
	\qquad P(w_n)=0.
	\]
	Moreover, since $\{y_n\}$ is bounded
	and $s_n\to0$, it follows from
	\eqref{eq:dilation-operator-bound}
	that $\{w_n\}$ is bounded in $X$.
	Finally, for any $k\in T_{w_n}\mathcal S$,
	the mass-preserving dilation implies that
	$(-s_n)\star k\in T_{y_n}\mathcal S$.
	Hence, by \eqref{eq:maximum-derivative}
	and \eqref{eq:dilation-operator-bound},
	\[
		|J'(w_n)[k]|
		=|dF(y_n)[(-s_n)\star k]|
		\leq
		\|dF(y_n)\|_{(T_{y_n}\mathcal S)^*}
		\|(-s_n)\star k\|_X
		\leq\varepsilon_n e^{|s_n|}\|k\|_X.
	\]
	Consequently,
	\[
	\|d(J|_{\mathcal S})(w_n)\|
	\leq\varepsilon_n e^{|s_n|}
	\to0,
	\]
	which completes the proof.
\end{proof}

\section{Compactness below the critical level}\label{sec:compactness}

In this section, we establish compactness results
below the critical energy level for use in the existence
proofs in Sections~\ref{sec:subcritical}
and~\ref{sec:critical-supercritical}.

\begin{lemma}[Annular regularity and Pohozaev identity]\label{lem:weak-pohozaev}
Let $w\in X$ be a nonzero weak solution of
\begin{equation}\label{eq:general-limiting-equation}
 -\operatorname{div}(|x|^{-2a}\nabla w)
 =\lambda|x|^{-2a}w+\beta|x|^{-bq}|w|^{q-2}w
 +\kappa|x|^{-bp}|w|^{p-2}w.
\end{equation}
Then $w\in L^\infty_{\rm loc}\cap H^2_{\rm loc}$ away from the
origin, and
\begin{equation}\label{eq:weak-pohozaev}
 P(w)=0,
 \qquad \lambda M(w)=(\delta_q-1)\beta B(w)<0.
\end{equation}
If $w\geq0$, then $w>0$ on $\mathbb R^N\setminus\{0\}$.
\end{lemma}

\begin{proof}
Let $\Omega'$ and $\Omega$ be bounded annuli
centered at the origin such that
$\Omega'\Subset\Omega\Subset\mathbb R^N\setminus\{0\}$.
Set
\begin{equation*}
	g(x)=|x|^{-2a},
	\qquad 2^*=\frac{2N}{N-2}.
\end{equation*}
There are constants $g_0,g_1>0$ such that
\begin{equation*}
 g_0\leq g(x)\leq g_1\qquad(x\in\Omega).
\end{equation*}
Since $d=1+a-b\in(0,1)$, we have
\begin{equation*}
 p=\frac{2N}{N-2d}<2^*,
 \qquad
 h:=\frac{2^*}{p-2}
 =\frac{N(N-2d)}{2d(N-2)}>\frac N2.
\end{equation*}
On $\Omega$, write \eqref{eq:general-limiting-equation} as
\begin{equation*}
 -\operatorname{div}(g\nabla w)=Vw,
 \qquad
 V=\lambda|x|^{-2a}
 +\beta|x|^{-bq}|w|^{q-2}
 +\kappa|x|^{-bp}|w|^{p-2}.
\end{equation*}
The embedding $H^1(\Omega)\hookrightarrow L^{2^*}(\Omega)$ gives
\begin{equation*}
 |w|^{p-2}\in L^h(\Omega),
 \qquad
 |w|^{q-2}\in L^{\frac{2^*}{q-2}}(\Omega)\subset L^h(\Omega).
\end{equation*}
Hence $V\in L^h(\Omega)$.

Let $s\geq 1$ and $L>0$. Choose
$\eta\in C_c^\infty(\Omega)$ such that
$0\leq\eta\leq1$ and $\eta\equiv1$ on $\Omega'$.
Set
\[
w_L=\min\{|w|,L\},
\qquad
f_L=\eta w w_L^{s-1}.
\]
Testing the weak equation with $\eta^2w w_L^{2s-2}$ gives
\[
\int g\nabla w\cdot\nabla
\bigl(\eta^2w w_L^{2s-2}\bigr)
=\int V\eta^2w^2w_L^{2s-2}.
\]
Expanding the left-hand side and applying the
Cauchy--Schwarz and Young inequalities, we obtain
\begin{equation}\label{eq:moser-caccioppoli}
	\|\nabla f_L\|_2^2
	\leq Cs^2\int |\nabla\eta|^2w^2w_L^{2s-2}
	+Cs^2\int |V|f_L^2,
\end{equation}
where $C>0$ is independent of $s$ and $L$.
Set
\[
r=\frac{2h}{h-1},
\qquad
\vartheta=\frac{N}{2h}\in(0,1).
\]
Since
\[
\frac1r
=\frac{\vartheta}{2^*}
+\frac{1-\vartheta}{2},
\]
the H\"older inequality and the interpolation inequality yield
\begin{equation*}
	\int |V|f_L^2
	\leq \|V\|_h\|f_L\|_r^2
	\leq C\|\nabla f_L\|_2^{2\vartheta}
	\|f_L\|_2^{2(1-\vartheta)}
	\leq \varepsilon\|\nabla f_L\|_2^2
	+C\varepsilon^{-\frac{N}{2h-N}}\|f_L\|_2^2
\end{equation*}
for every $\varepsilon>0$, where the last
inequality follows from the Young inequality.
Choosing $\varepsilon=\frac{1}{2C s^2}$
in \eqref{eq:moser-caccioppoli} and
absorbing the gradient term, we obtain
\begin{equation}\label{eq:moser-absorbed}
	\|\nabla f_L\|_2^2
	\leq Cs^\ell\int
	\bigl(|\nabla\eta|^2+\eta^2\bigr)
	w^2w_L^{2s-2},
\end{equation}
where
$
\ell=2+\frac{2N}{2h-N}.
$
Here and below, $C>0$ denotes a constant
independent of $s$ and $L$.
Assume that $w\in L^{2s}(\Omega)$.
Applying the Sobolev inequality to $f_L$
and then letting $L\to\infty$, we deduce
from \eqref{eq:moser-absorbed} and
Fatou's lemma that
\begin{equation}\label{eq:moser-step}
	\|w\|_{L^{2^*s}(\Omega')}
	\leq
	\bigl(Cs^\ell
	(1+\|\nabla\eta\|_\infty^2)\bigr)^{\frac1{2s}}
	\|w\|_{L^{2s}(\Omega)}.
\end{equation}

Choose a decreasing sequence of annuli
$\{\Omega_j\}_{j\geq0}$ such that
\[
\Omega'\Subset\Omega_{j+1}\Subset\Omega_j\Subset\Omega,
\qquad
\bigcap_{j=0}^\infty\Omega_j=\overline{\Omega'}.
\]
For each $j\geq0$, let
$\eta_j\in C_c^\infty(\Omega_j)$ satisfy
\[
0\leq\eta_j\leq1,\qquad
\eta_j\equiv1\ \text{on }\Omega_{j+1},
\qquad
\|\nabla\eta_j\|_\infty\leq C2^j.
\]
Set
\[
\sigma=\frac{N}{N-2}>1,
\qquad s_j=\sigma^j.
\]
Since $2^*s_j=2s_{j+1}$,
applying \eqref{eq:moser-step} with
$s=s_j$ and $\eta=\eta_j$ yields
\[
\|w\|_{L^{2s_{j+1}}(\Omega_{j+1})}
\leq
\bigl(C4^js_j^\ell\bigr)^{1/(2s_j)}
\|w\|_{L^{2s_j}(\Omega_j)}.
\]
Moreover,
\[
\sum_{j=0}^\infty
\frac{\log C+j\log4+\ell\log s_j}{2s_j}
<\infty,
\]
since $s_j=\sigma^j$ with $\sigma>1$.
Consequently, iterating the preceding
inequality and letting $j\to\infty$,
we obtain
\[
\|w\|_{L^\infty(\Omega')}
\leq C\|w\|_{L^2(\Omega)}.
\]
As $\Omega'$ is arbitrary, it follows that
\[
w\in L^\infty_{\mathrm{loc}}
(\mathbb R^N\setminus\{0\}).
\]
Dividing the equation by $g(x)=|x|^{-2a}$ on $\Omega$ gives
\[
-\Delta w
+2a\frac{x}{|x|^2}\cdot\nabla w
=\lambda w
+\beta|x|^{2a-bq}|w|^{q-2}w
+\kappa|x|^{2a-bp}|w|^{p-2}w.
\]
Its right-hand side is locally bounded. Interior estimates give
\begin{equation*}
 w\in H^2_{\rm loc}(\mathbb R^N\setminus\{0\})
 \cap W^{2,r}_{\rm loc}(\mathbb R^N\setminus\{0\})
 \qquad(1<r<\infty);
\end{equation*}
see \cite[Chapters 8--9]{GilbargTrudinger2001}.

Let $0<2\varepsilon<R$. Choose a cutoff function
$\chi_{\varepsilon,R}\in C_c^\infty(\mathbb R^N\setminus\{0\})$
satisfying $0\leq\chi_{\varepsilon,R}\leq1$ and
\begin{equation*}
	\chi_{\varepsilon,R}=
	\begin{cases}
		1, & \text{if } \varepsilon\leq |x|\leq R,\\
		0, & \text{if } |x|\leq\varepsilon/2
		\text{ or } |x|\geq2R,
	\end{cases}
	\qquad
	|x|\,|\nabla\chi_{\varepsilon,R}|\leq C.
\end{equation*}
The function
\[
\varphi_{\varepsilon,R}=\chi_{\varepsilon,R}x\cdot\nabla w
\]
belongs to $H_0^1(B_{2R}\setminus B_{\varepsilon/2})$ and is an admissible test function for the weak formulation.
Put $\chi=\chi_{\varepsilon,R}$ and
\begin{equation*}
 \mathcal A_{\varepsilon,R}
 =\{\varepsilon/2<|x|<\varepsilon\}
 \cup\{R<|x|<2R\}.
\end{equation*}
Integration by parts gives
\begin{align*}
 \int |x|^{-2a}\nabla w\cdot\nabla(\chi x\cdot\nabla w)
 &=-\frac{N-2-2a}{2}
 \int\chi|x|^{-2a}|\nabla w|^2+\mathcal R_A,\\
 \lambda\int\chi|x|^{-2a}w(x\cdot\nabla w)
 &=-\frac{N-2a}{2}\lambda
 \int\chi|x|^{-2a}w^2+\mathcal R_M,\\
 \beta\int\chi|x|^{-bq}|w|^{q-2}w(x\cdot\nabla w)
 &=-\frac{N-bq}{q}\beta
 \int\chi|x|^{-bq}|w|^q+\mathcal R_B,\\
 \kappa\int\chi|x|^{-bp}|w|^{p-2}w(x\cdot\nabla w)
 &=-\frac{N-bp}{p}\kappa
 \int\chi|x|^{-bp}|w|^p+\mathcal R_C,
\end{align*}
where
\begin{align*}
 \mathcal R_A
 &=\int |x|^{-2a}(x\cdot\nabla w)
 (\nabla\chi\cdot\nabla w)
 -\frac12\int |x|^{-2a}(x\cdot\nabla\chi)|\nabla w|^2,\\
 \mathcal R_M
 &=-\frac{\lambda}{2}\int
 |x|^{-2a}(x\cdot\nabla\chi)w^2,\\
 \mathcal R_B
 &=-\frac{\beta}{q}\int
 |x|^{-bq}(x\cdot\nabla\chi)|w|^q,\\
 \mathcal R_C
 &=-\frac{\kappa}{p}\int
 |x|^{-bp}(x\cdot\nabla\chi)|w|^p.
\end{align*}
Set
\begin{equation*}
 e_w=|x|^{-2a}|\nabla w|^2
 +|\lambda||x|^{-2a}w^2
 +\beta|x|^{-bq}|w|^q
 +\kappa|x|^{-bp}|w|^p.
\end{equation*}
Then
\begin{equation*}
 |\mathcal R_A|+|\mathcal R_M|+|\mathcal R_B|+|\mathcal R_C|
 \leq C\int_{\mathcal A_{\varepsilon,R}}e_w.
\end{equation*}
Lemma \ref{lem:weighted-gn} and \eqref{eq:sharp-ckn} imply
$e_w\in L^1(\mathbb R^N)$. Therefore
\begin{equation*}
 \lim_{\varepsilon\to0^+}\int_{|x|<\varepsilon}e_w=0,
 \qquad
 \lim_{R\to+\infty}\int_{|x|>R}e_w=0.
\end{equation*}
Letting $\varepsilon\to0^+$ and $R\to+\infty$ gives
\begin{equation}\label{eq:spatial-pohozaev-cutoff}
 \frac{N-2-2a}{2}A(w)
 =\frac{N-2a}{2}\lambda M(w)
 +\frac{N-bq}{q}\beta B(w)
 +\frac{N-bp}{p}\kappa C(w).
\end{equation}
Testing \eqref{eq:general-limiting-equation} with $w$ yields
\begin{equation}\label{eq:weak-nehari-identity}
	A(w)=\lambda M(w)+\beta B(w)+\kappa C(w).
\end{equation}
Set
\[
c_a=\frac{N-2-2a}{2}.
\]
By the definitions of $p$ and $\delta_q$, we have
\[
\frac{N-2a}{2}=c_a+1,
\qquad
\frac{N-bq}{q}=c_a+1-\delta_q,
\qquad
\frac{N-bp}{p}=c_a.
\]
Subtracting $c_a$ times
\eqref{eq:weak-nehari-identity} from
\eqref{eq:spatial-pohozaev-cutoff}, we obtain
\begin{equation}\label{eq:weak-multiplier-identity}
	\lambda M(w)=(\delta_q-1)\beta B(w).
\end{equation}
Consequently,
\[
	P(w)
	=A(w)-\delta_q\beta B(w)-\kappa C(w)
	=\lambda M(w)+(1-\delta_q)\beta B(w)
	=0.
\]
The function
\[
r\longmapsto\delta_r
=\frac{N-2a+2b}{2}-\frac Nr
\]
is strictly increasing on $(0,\infty)$, so $q<p$ implies
$\delta_q<\delta_p=1$.
Since $w\not\equiv0$, both $M(w)$ and $B(w)$ are positive,
so \eqref{eq:weak-multiplier-identity} yields $\lambda<0$.

Finally, assume $w\geq0$. On every annulus,
\begin{equation*}
 -\operatorname{div}(|x|^{-2a}\nabla w)
 -\lambda|x|^{-2a}w
 =\beta|x|^{-bq}w^{q-1}
 +\kappa|x|^{-bp}w^{p-1}\geq0.
\end{equation*}
The operator is uniformly elliptic there, and its zero-order
coefficient is positive. If $w(x_0)=0$ for some $x_0\ne0$, choose
a bounded connected domain compactly contained in
$\mathbb R^N\setminus\{0\}$ that contains $x_0$ and a point where
$w>0$. The strong maximum principle gives a contradiction.
Thus $w>0$ on $\mathbb R^N\setminus\{0\}$.
\end{proof}

\begin{lemma}\label{lem:scalar-splitting}
Let $w_n\rightharpoonup w$ in $X$ and set $z_n=w_n-w$.  Then, after passing to a subsequence,
\begin{equation*}
 A(w_n)=A(w)+A(z_n)+o(1),
\end{equation*}
\begin{equation*}
 M(w_n)=M(w)+M(z_n)+o(1),
\end{equation*}
\begin{equation}\label{eq:B-compact}
 B(w_n)=B(w)+o(1),
\end{equation}
and
\begin{equation}\label{eq:C-splitting}
 C(w_n)=C(w)+C(z_n)+o(1).
\end{equation}
\end{lemma}

\begin{proof}
The first two identities follow from weak convergence
in the two Hilbert structures.
By Proposition~\ref{prop:compact-embedding},
$w_n\to w$ strongly in the weighted $L^q$ space,
which yields \eqref{eq:B-compact}.
Applying the Brezis--Lieb lemma \cite{BrezisLieb1983}
to $|x|^{-b}w_n$ gives \eqref{eq:C-splitting}.
\end{proof}

The following proposition establishes mass recovery for nontrivial weak limits of bounded constrained Palais--Smale sequences and quantifies the possible energy loss due to critical concentration.

\begin{proposition}[Mass recovery and critical energy defect]
	\label{prop:compactness-alternative}
	Let $\{w_n\}\subset\mathcal S$ be a bounded
	Palais--Smale sequence for $J|_{\mathcal S}$
	at level $c$, satisfying $P(w_n)\to0$.
	Then, up to a subsequence, there exist
	$w\in X$ and $\lambda\in\mathbb R$ such that
	$
	w_n\rightharpoonup w\ \text{in }X,
	\lambda_n\to\lambda,
	$
	where $\lambda_n$ denotes the associated
	Lagrange multiplier.
	
	If $w\not\equiv0$, then $w$ satisfies
	\eqref{eq:general-limiting-equation}, and
	\begin{equation}\label{eq:early-mass-recovery}
		M(w)=1,\qquad M(w_n-w)\to0,
		\qquad
		\lambda=(\delta_q-1)\beta B(w)<0.
	\end{equation}
	Whether or not $w\equiv0$, there exists
	$\Lambda\geq0$ such that
	\begin{equation}\label{eq:critical-defect-splitting}
		A(w_n-w)\to\Lambda,\qquad
		\kappa C(w_n-w)\to\Lambda,
		\qquad
		c=J(w)+\frac{d}{N}\Lambda.
	\end{equation}
	Furthermore, either $\Lambda=0$, or
	\begin{equation}\label{eq:critical-defect-lower}
		\Lambda\geq
		S(a,b)^{\frac{N}{2d}}
		\kappa^{-\frac{N-2d}{2d}},
		\qquad
		J(w)\leq c-T_\kappa.
	\end{equation}
	In particular, if $\Lambda=0$ and $c\neq0$,
	then $w_n\to w$ strongly in $X$ and
	$w\in\mathcal S$.
\end{proposition}

\begin{proof}
	Since $\{w_n\}$ is bounded in $X$, passing to a
	subsequence if necessary, we may assume that
	$
	w_n\rightharpoonup w\ \text{in }X.
	$
	By the constrained Palais--Smale condition and
	the standard Lagrange multiplier principle,
	there exists a sequence $\{\lambda_n\}\subset\mathbb R$
	such that, for every $\varphi\in X$,
	\begin{align}
		&\int |x|^{-2a}\nabla w_n\cdot\nabla\varphi
		-\lambda_n\int |x|^{-2a}w_n\varphi
		-\beta\int |x|^{-bq}|w_n|^{q-2}w_n\varphi
		\notag\\
		&\qquad
		-\kappa\int |x|^{-bp}|w_n|^{p-2}w_n\varphi
		=o(1)\|\varphi\|_X.
		\label{eq:PS-equation}
	\end{align}
    Testing \eqref{eq:PS-equation} with $w_n$ and using $M(w_n)=1$ yields
	\begin{equation}\label{eq:sequence-multiplier}
		\lambda_n
		=P(w_n)+(\delta_q-1)\beta B(w_n)+o(1)
		\longrightarrow
	    (\delta_q-1)\beta B(w)=:\lambda,
	\end{equation}
	where the convergence follows from
	$P(w_n)\to0$ and
	Proposition~\ref{prop:compact-embedding}.
	
	Since $\{|x|^{-b}w_n\}$ is bounded in $L^p$ and converges almost everywhere to $|x|^{-b}w$, the corresponding critical nonlinearities converge weakly in $X^*$. Using this convergence, the strong convergence of the subcritical nonlinearities in $X^*$, and $\lambda_n\to\lambda$, we pass to the limit in \eqref{eq:PS-equation} and conclude that $w$ satisfies \eqref{eq:general-limiting-equation}.

	If $w\neq0$, Lemma~\ref{lem:weak-pohozaev} yields
\begin{equation*}
	P(w)=0,\qquad \lambda M(w)=(\delta_q-1)\beta B(w)<0.
\end{equation*}
Comparing these identities with \eqref{eq:sequence-multiplier},
we obtain $\lambda<0$ and $M(w)=1$. The Hilbert space splitting
of the mass then gives \eqref{eq:early-mass-recovery}.
If $w=0$, the identity $P(w)=0$ also holds.
	
Set $z_n=w_n-w$. Lemma~\ref{lem:scalar-splitting},
together with $P(w_n)\to0$ and $P(w)=0$, gives
\[
A(z_n)-\kappa C(z_n)=o(1),
\qquad
J(w_n)=J(w)+\frac12A(z_n)
-\frac{\kappa}{p}C(z_n)+o(1).
\]
Since $\{z_n\}$ is bounded in $X$, passing to a
subsequence, we may assume that
\[
A(z_n)\to\Lambda\geq0,
\qquad \kappa C(z_n)\to\Lambda.
\]
Since $J(w_n)\to c$, it follows that
\[
c=J(w)+\frac dN\Lambda,
\]
which proves \eqref{eq:critical-defect-splitting}.
By \eqref{eq:sharp-ckn},
$
S(a,b)\left(\frac{\Lambda}{\kappa}\right)^{2/p}
\leq\Lambda.
$
Since $p>2$, it follows that either $\Lambda=0$ or
$
\Lambda
\geq S(a,b)^{\frac{N}{2d}}
\kappa^{-\frac{N-2d}{2d}}.
$
In the latter case, \eqref{eq:critical-defect-splitting} gives
\[
J(w)=c-\frac{d}{N}\Lambda\leq c-T_\kappa,
\]
which proves \eqref{eq:critical-defect-lower}.
Finally, if $\Lambda=0$ and $w=0$, the energy splitting
gives $c=0$. Hence, when $\Lambda=0$ and $c\neq0$,
we have $w\neq0$. The vanishing gradient defect, together
with \eqref{eq:early-mass-recovery}, then yields strong
convergence in $X$.
\end{proof}

\begin{lemma}[Mass-dependent coercivity]\label{lem:negative-lower-bound}
Let $2<q<q_c$. Set
\begin{equation*}
 a_0=\frac dN,\qquad
 b_0=\left(\frac1q-\frac{\delta_q}{p}\right)K_q,
 \qquad s_0=\frac\theta2,
 \qquad C_0=a_0\frac{1-s_0}{s_0}
 \left(\frac{b_0s_0}{a_0}\right)^{\frac{1}{1-s_0}}.
\end{equation*}
Every $w\in X$ with $P(w)=0$ satisfies
\begin{equation}\label{eq:mass-dependent-coercivity}
 J(w)\geq a_0A(w)-b_0\beta M(w)^{q(1-\delta_q)/2}A(w)^{s_0}
 \geq-C_0\beta^{2/(2-\theta)}M(w)^{q(1-\delta_q)/(2-\theta)}.
\end{equation}
In particular, $J$ is coercive as a function of $A$ on
$\{w:P(w)=0, M(w)\leq1\}$.
\end{lemma}

\begin{proof}
Combining \eqref{eq:energy-on-P-A} with
\eqref{eq:weighted-gn}, we obtain
\[
J(w)\geq a_0A(w)
-b_0\beta M(w)^{q(1-\delta_q)/2}A(w)^{s_0}.
\]

Since $0<s_0<1$, minimizing the right-hand side over $A(w)\geq0$
yields the second estimate in \eqref{eq:mass-dependent-coercivity}.
The first estimate also gives coercivity for $M(w)\leq1$.
\end{proof}

\section{Two solutions in the mass-subcritical case}\label{sec:subcritical}

In this section, we prove the existence of two
solutions with opposite energies in the
mass-subcritical regime, corresponding to a local
minimizer and an upper Pohozaev branch minimizer.

Throughout this section, we assume that
\begin{equation*}
	2<q<q_c,
	\qquad
	0<\beta<\beta_0.
\end{equation*}
The inequalities \eqref{eq:sharp-ckn} and \eqref{eq:weighted-gn} imply
\begin{equation}\label{eq:energy-basic-lower}
	J(w)\geq f_{\beta}(A(w)^{1/2}),
	\qquad w\in\mathcal S.
\end{equation}
Since $\beta<\beta_0\leq\beta_{\rm b}$,
the definitions of $R_*$ and $\beta_{\rm b}$ yield
\begin{equation}\label{eq:barrier-positive}
	f_{\beta}(R_*)\geq\frac14R_*^2>0.
\end{equation}
Define
\begin{equation*}
	\mathcal G
	=\{w\in\mathcal S:A(w)^{1/2}<R_*\},
	\qquad
	m^+=\inf_{\mathcal G}J.
\end{equation*}

\begin{proposition}\label{prop:local-minimizer}
Every minimizing sequence for $m^+$ has a strongly convergent
subsequence in $X$. This level is attained by a nonnegative
local minimizer $w^+$ of $J|_{\mathcal S}$, and
\begin{equation}\label{eq:local-min-energy}
 J(w^+)=m^+<0,
 \qquad w^+\in\mathcal P^+,
 \qquad m^+=\inf_{\mathcal P^+}J.
\end{equation}
Every point of $\mathcal P^+$ belongs to $\mathcal G$.
\end{proposition}

\begin{proof}
By \eqref{eq:energy-basic-lower} and the definition of
$\mathcal G$, we have
$
m^+\geq\min_{0\leq s\leq R_*}f_\beta(s)>-\infty.
$
For any $w\in\mathcal S$, as $t\to-\infty$,
\[
J(t\star w)
=-\frac{\beta}{q}e^{\theta t}B(w)
+o(e^{\theta t})<0,
\qquad
A(t\star w)^{1/2}\to0.
\]
Hence $m^+<0$.
By the continuity of $f_\beta$ and
\eqref{eq:barrier-positive}, there exists
$\varepsilon_*>0$ such that
\[
J(w)>0
\quad\text{whenever}\quad
R_*-\varepsilon_*
\leq A(w)^{1/2}\leq R_*,
\qquad w\in\mathcal S.
\]
Thus, any minimizing sequence
$\{w_n\}\subset\mathcal G$ satisfies
$
A(w_n)^{1/2}<R_*-\varepsilon_*
$
for all sufficiently large $n$.
Since $\{w_n\}$ is bounded in $X$,
up to a subsequence, we have
$
w_n\rightharpoonup w\ \text{in }X,
A(w)^{1/2}\leq R_*-\varepsilon_*.
$
Set $z_n=w_n-w$. Then
\[
A(z_n)=A(w_n)-A(w)+o(1)
\leq R_*^2+o(1).
\]
It follows from Lemma~\ref{lem:scalar-splitting}
and \eqref{eq:sharp-ckn} that
\begin{equation}\label{eq:direct-well-defect}
	\begin{aligned}
		J(w_n)-J(w)
		&=\frac12 A(z_n)-\frac{\kappa}{p}C(z_n)+o(1)\\
		&\geq
		\left(
		\frac12-\frac{\kappa R_*^{p-2}}
		{pS(a,b)^{p/2}}
		\right)A(z_n)+o(1)\\
		&=c_2A(z_n)+o(1),
	\end{aligned}
\end{equation}
where $c_2>0$ by \eqref{eq:R-star}.
Consequently, $J(w)\leq m^+<0$, and hence
$w\neq0$.

Write $M(w)=\sigma^2\leq1$. If $\sigma<1$, define
\begin{equation*}
	\widehat w(x)=\sigma^{(N-2-2a)/2}w(\sigma x).
\end{equation*}
The critical scaling properties yield
\begin{equation*}
	M(\widehat w)=1,\quad A(\widehat w)=A(w),\quad
	C(\widehat w)=C(w),\quad B(\widehat w)=\sigma^{-q(1-\delta_q)}B(w).
\end{equation*}
Since $q(1-\delta_q)>0$ and $B(w)>0$,
we have $\widehat w\in\mathcal G$ and
\[
J(\widehat w)
=J(w)-\frac{\beta}{q}
\bigl(\sigma^{-q(1-\delta_q)}-1\bigr)B(w)
<J(w)\leq m^+,
\]
a contradiction. Thus $M(w)=1$ and $J(w)=m^+$. Together with
$J(w_n)\to m^+$, \eqref{eq:direct-well-defect}
and the mass splitting yield $A(z_n)\to0$ and
$M(z_n)\to0$, respectively. Hence $w_n\to w$
strongly in $X$.
As an interior minimizer of $J$ on $\mathcal G$
relative to $\mathcal S$, $w$ is a constrained
critical point of $J$ on $\mathcal S$.
Lemma~\ref{lem:weak-pohozaev} gives $P(w)=0$,
and $J(w)=m^+<0$ then implies $w\in\mathcal P^+$
by Lemma~\ref{lem:subcritical-fibers}.
Since replacing $w$ by $|w|$ preserves
$A$, $M$, $B$, and $C$, the minimizer may be
chosen nonnegative.

Finally, let $z\in\mathcal P^+$. Since the lower
fiber critical point occurs at $t=0$, we have
$J(t\star z)<0$ for all $t\leq0$.
On the other hand, \eqref{eq:fiber-positive-at-R}
gives $J(t_R\star z)>0$, where
$
e^{t_R}A(z)^{1/2}=R_*.
$
Hence $t_R>0$, which implies
$A(z)^{1/2}<R_*$ and thus
$\mathcal P^+\subset\mathcal G$.
Together with the minimizer obtained above,
this yields the last identity in
\eqref{eq:local-min-energy}.
\end{proof}

We define the upper-branch energy level by
\begin{equation*}
	m^-:=\inf_{w\in\mathcal P^-}J(w).
\end{equation*}

\begin{lemma}\label{lem:upper-positive}
	The upper-branch energy level satisfies
	\[
	m^-\geq \frac{R_*^2}{4}>0.
	\]
\end{lemma}

\begin{proof}
For $w\in\mathcal P^-$, Lemma \ref{lem:subcritical-fibers} and
\eqref{eq:fiber-positive-at-R} give
\begin{equation*}
 J(w)=\max_tJ(t\star w)\geq f_\beta(R_*)\geq R_*^2/4.
\end{equation*}
Take the infimum over $\mathcal P^-$.
\end{proof}

For $\varepsilon>0$, define the rescaling of the finite-mass
extremal $W$ by
\begin{equation*}
	W_{\varepsilon}(x)
	=\varepsilon^{-(N-2-2a)/2}W(x/\varepsilon).
\end{equation*}
Then
\begin{equation*}
 A(W_{\varepsilon})=A(W),
 \quad C(W_{\varepsilon})=C(W),
 \quad M(W_{\varepsilon})=\varepsilon^2M(W),
 \quad B(W_{\varepsilon})=\varepsilon^{q(1-\delta_q)}B(W).
\end{equation*}
Let
\begin{equation*}
 J_0(w)=\frac12A(w)-\frac{\kappa}{p}C(w).
\end{equation*}
The function $W$ is an extremal for the sharp CKN
inequality \eqref{eq:sharp-ckn}, so
\begin{equation}\label{eq:pure-critical-maximum}
\max_{s>0}J_0(sW)=T_\kappa.
\end{equation}

\begin{lemma}\label{lem:path-estimate}
	Let $w^+$ be the minimizer obtained in
	Proposition~\ref{prop:local-minimizer}.
	Then the upper Pohozaev branch level satisfies
	\[
	m^-<m^++T_\kappa.
	\]
\end{lemma}

\begin{proof}
Recall that the Lagrange multiplier associated with
$w^+$ satisfies
$
\lambda=(\delta_q-1)\beta B(w^+)<0.
$
For $\varepsilon>0$ and $s\geq0$, set
\[
h_{\varepsilon,s}=w^++sW_\varepsilon,
\qquad
\tau_{\varepsilon,s}=M(h_{\varepsilon,s})^{1/2},
\]
and define the normalized path
\[
v_{\varepsilon,s}(x)
=\tau_{\varepsilon,s}^{(N-2-2a)/2}
h_{\varepsilon,s}(\tau_{\varepsilon,s}x).
\]
The scaling properties of the critical dilation yield
\begin{equation}\label{eq:path-normalization-identities}
	M(v_{\varepsilon,s})=1,
	\quad A(v_{\varepsilon,s})=A(h_{\varepsilon,s}),
	\quad C(v_{\varepsilon,s})=C(h_{\varepsilon,s}),
	\quad B(v_{\varepsilon,s})=\tau_{\varepsilon,s}^{-q(1-\delta_q)}B(h_{\varepsilon,s}).
\end{equation}
Set
\begin{equation*}
	m_{\varepsilon}=\int |x|^{-2a}w^+W_{\varepsilon},
	\qquad
	Q_{\varepsilon}=\int |x|^{-bq}(w^+)^{q-1}W_{\varepsilon}.
\end{equation*}
The H\"older inequality and the scaling properties of $W_\varepsilon$ give
\[
0\leq m_\varepsilon
\leq\varepsilon M(W)^{1/2},
\qquad
0\leq Q_\varepsilon
\leq\varepsilon^{1-\delta_q}
B(w^+)^{(q-1)/q}B(W)^{1/q}.
\]
Moreover, expanding $M(h_{\varepsilon,s})$ gives
\begin{equation}\label{eq:tau-square}
	\tau_{\varepsilon,s}^2
	=1+2sm_\varepsilon+s^2\varepsilon^2M(W).
\end{equation}

For $r\geq2$ and $x,y\geq0$,
\begin{equation}\label{eq:convex-power}
	(x+y)^r\geq x^r+y^r+rx^{r-1}y.
\end{equation}
Testing the equation for $w^+$ with $W_\varepsilon$
gives
\begin{equation}\label{eq:cross-equation}
	\int |x|^{-2a}\nabla w^+\cdot\nabla W_\varepsilon
	=\lambda m_\varepsilon+\beta Q_\varepsilon
	+\kappa\int |x|^{-bp}
	(w^+)^{p-1}W_\varepsilon.
\end{equation}
Applying \eqref{eq:convex-power} to the $q$- and
$p$-power terms and using \eqref{eq:cross-equation},
we obtain, uniformly for $s$ in bounded intervals,
\begin{align}
	J(v_{\varepsilon,s})-J(w^+)-J_0(sW)
	&\leq s\lambda m_\varepsilon
	+\frac{\beta}{q}
	(1-\tau_{\varepsilon,s}^{-q(1-\delta_q)})
	B(w^+)\notag\\
	&\quad+\beta s
	(1-\tau_{\varepsilon,s}^{-q(1-\delta_q)})
	Q_\varepsilon-\frac{\beta}{q}
	\tau_{\varepsilon,s}^{-q(1-\delta_q)}
	s^q\varepsilon^{q(1-\delta_q)}B(W).
	\label{eq:path-expansion}
\end{align}
To estimate the first two terms, we use
\[
1-x^{-q(1-\delta_q)/2}
\leq\frac{q(1-\delta_q)}{2}(x-1),
\qquad x\geq1.
\]
Since $\tau_{\varepsilon,s}^2\geq1$ and
\[
-\lambda=\beta(1-\delta_q)B(w^+),
\]
it follows from \eqref{eq:tau-square} that
\begin{equation}\label{eq:first-order-cancellation}
	s\lambda m_\varepsilon
	+\frac{\beta}{q}
	(1-\tau_{\varepsilon,s}^{-q(1-\delta_q)})
	B(w^+)\leq
	-\frac{\lambda}{2}s^2\varepsilon^2M(W).
\end{equation}
Moreover, the estimates for $m_\varepsilon$ and
$Q_\varepsilon$, together with \eqref{eq:tau-square},
yield
\begin{equation}\label{eq:Q-error}
	0\leq
	\beta s(1-\tau_{\varepsilon,s}^{-q(1-\delta_q)})
	Q_\varepsilon
	\leq C\left(
	\varepsilon^{2-\delta_q}
	+\varepsilon^{3-\delta_q}
	\right),
\end{equation}
uniformly for $s$ in bounded intervals.
We next establish an estimate that holds uniformly down to $s=0$. Since $2<q<q_c<p$, we have
\[
0<q(1-\delta_q)
=N-\frac{N-2d}{2}q<2d<2.
\]
Moreover,
\[
	2-\delta_q-q(1-\delta_q)
	=1-d+(q-2)
	\left(\frac{N(q-1)}{2q}-d\right)>1-d>0,
\]
since $\frac{N(q-1)}{2q}>\frac N4>d$.
Consequently,
\[
\zeta_0:=
\min\left\{
2-q(1-\delta_q),
2-\delta_q-q(1-\delta_q)
\right\}>0.
\]
For any fixed $L>0$, \eqref{eq:tau-square} and
$m_\varepsilon=O(\varepsilon)$ yield
\[
1\leq\tau_{\varepsilon,s}^2
\leq1+C_L\varepsilon,
\qquad 0\leq s\leq L,
\]
where $C_L>0$ is independent of $s$ and
$\varepsilon$. Consequently,
$
\tau_{\varepsilon,s}^{-q(1-\delta_q)}
\geq\frac12
$
for sufficiently small $\varepsilon$,
uniformly in $s\in[0,L]$.
Combining \eqref{eq:path-expansion},
\eqref{eq:first-order-cancellation}, and
\eqref{eq:Q-error}, we obtain
\begin{equation}\label{eq:path-final-expansion}
	J(v_{\varepsilon,s})
	\leq m^++J_0(sW)
	-\frac{\beta}{2q}s^qB(W)
	\varepsilon^{q(1-\delta_q)}+C_L\varepsilon^{q(1-\delta_q)+\zeta_0},
	\quad 0\leq s\leq L.
\end{equation}

Since \eqref{eq:path-final-expansion} is valid only for
bounded $s$ and its negative correction degenerates as
$s\to0$, we first choose the path parameters independently
of $\varepsilon$ and then treat different ranges of $s$
separately. By
\eqref{eq:path-normalization-identities}
and the nonnegativity of $w^+$ and $W_\varepsilon$, we have
\[
J(v_{\varepsilon,s})
\leq H(s):=
A(w^+)+s^2A(W)-\frac{\kappa}{p}s^pC(W),
\qquad s\geq0.
\]
Let $s_*>0$ be the unique maximizer in
\eqref{eq:pure-critical-maximum}.
Since $H(s)\to-\infty$ as $s\to\infty$,
we may fix $s_1>s_*$ sufficiently large
such that
\[
H(s)\leq m^++T_\kappa-1,
\qquad s\geq2s_1.
\]
Choose $s_0\in(0,s_*)$ and set
\[
\zeta=\frac14\min\left\{
1,\,
T_\kappa-J_0(s_0W),\,
T_\kappa-J_0(s_1W)
\right\}>0.
\]
By the monotonicity of the pure critical fiber,
\begin{equation}\label{eq:pure-gap}
	J_0(sW)\leq T_\kappa-3\zeta,
	\qquad
	s\in[0,s_0]\cup[s_1,2s_1].
\end{equation}
Applying \eqref{eq:path-final-expansion}
with $L=2s_1$, we choose
$\varepsilon_0\in(0,1)$ sufficiently small
such that the estimate holds for
$0<\varepsilon\leq\varepsilon_0$ and
\[
C_{2s_1}\varepsilon_0^{\zeta_0}
\leq\frac{\beta}{4q}s_0^qB(W),
\qquad
C_{2s_1}
\varepsilon_0^{q(1-\delta_q)+\zeta_0}
\leq\zeta.
\]
For $s\in[s_0,s_1]$,
\[
J(v_{\varepsilon,s})
\leq m^++T_\kappa
-\frac{\beta}{4q}s_0^qB(W)
\varepsilon^{q(1-\delta_q)}.
\]
On $[0,s_0]\cup[s_1,2s_1]$,
\eqref{eq:path-final-expansion} and
\eqref{eq:pure-gap} give
\[
J(v_{\varepsilon,s})
\leq m^++T_\kappa-2\zeta.
\]
For $s\geq2s_1$, the bound involving $H$
yields
\[
J(v_{\varepsilon,s})
\leq m^++T_\kappa-1.
\]
Therefore, for every
$0<\varepsilon\leq\varepsilon_0$,
\begin{equation}\label{eq:path-strict-upper}
	\sup_{s\geq0}J(v_{\varepsilon,s})
	<m^++T_\kappa.
\end{equation}

It remains to show that the path intersects
$\mathcal P^-$. Let
$t_2(v_{\varepsilon,s})$ be the upper fiber
parameter given by
Lemma~\ref{lem:subcritical-fibers}.
Since $v_{\varepsilon,0}=w^+\in\mathcal P^+$,
we have
$
t_2(v_{\varepsilon,0})>0.
$
The Pohozaev identity at the upper fiber
critical point gives
\begin{equation}\label{eq:upper-root-ratio}
	e^{(p-2)t_2(v_{\varepsilon,s})}
	\leq
	\frac{A(v_{\varepsilon,s})}
	{\kappa C(v_{\varepsilon,s})}.
\end{equation}
Moreover,
\[
A(v_{\varepsilon,s})
\leq2A(w^+)+2s^2A(W),
\qquad
C(v_{\varepsilon,s})\geq s^pC(W).
\]
Consequently, for $s\geq1$,
\[
e^{(p-2)t_2(v_{\varepsilon,s})}
\leq Cs^{2-p},
\]
where $C>0$ is independent of
$\varepsilon$ and $s$. Thus
\[
t_2(v_{\varepsilon,s})
\leq\frac{\log C}{p-2}-\log s.
\]
In particular, there exists $L_1>0$,
independent of $\varepsilon$, such that
$
t_2(v_{\varepsilon,L_1})<0.
$
By continuity, there exists $s_{\varepsilon}\in(0,L_1)$
such that
\begin{equation*}
	t_2(v_{\varepsilon,s_{\varepsilon}})=0,
	\qquad
	v_{\varepsilon,s_{\varepsilon}}\in\mathcal P^-.
\end{equation*}
Combining this with \eqref{eq:path-strict-upper},
we conclude that
\[
m^-
\leq J(v_{\varepsilon,s_\varepsilon})
\leq\sup_{s\geq0}J(v_{\varepsilon,s})
<m^++T_\kappa.
\]
\end{proof}

\begin{proposition}\label{prop:upper-solution}
The level $m^-$ is attained by a nonnegative constrained critical
point $w^-\in\mathcal P^-$ of $J|_{\mathcal S}$.
\end{proposition}
\begin{proof}
Let $\{w_n\}\subset\mathcal P^-$ be a minimizing sequence
for $m^-$. By \eqref{eq:mass-dependent-coercivity}, the
sequence $\{w_n\}$ is bounded in $X$. Proposition
\ref{prop:full-PS} then yields a bounded Palais--Smale
sequence, still denoted by $\{w_n\}$, satisfying
\eqref{eq:full-PS-sequence} with $c=m^-$.
	
If critical loss occurs, Proposition \ref{prop:compactness-alternative}
and Lemma \ref{lem:path-estimate} give
\begin{equation*}
	J(w)\leq m^--T_\kappa<m^+<0.
\end{equation*}
Thus $w\neq0$. By
\eqref{eq:early-mass-recovery}, we have $M(w)=1$, while
Lemma~\ref{lem:weak-pohozaev} gives $P(w)=0$.
Since $J(w)<0$, the fiber geometry in Lemma
\ref{lem:subcritical-fibers} implies that
$w\in\mathcal P^+$. Therefore, Proposition
\ref{prop:local-minimizer} yields
\[
m^+\leq J(w)<m^+,
\]
which is impossible.
Since $m^->0$, the compactness alternative yields, up to a subsequence, $w_n\to w^-$ strongly in $X$ for a nontrivial limit $w^-$.
Moreover,
$
P(w^-)=0,
$
and \eqref{eq:uniform-fiber-gap} implies that
$w^-\in\mathcal P^-$. Consequently,
\[
J(w^-)=m^-.
\]
Finally, replacing $w^-$ by $|w^-|$ preserves the four
scalar quantities appearing in $J$ and $P$. Hence
$|w^-|$ is also a minimizer on $\mathcal P^-$. By
Lemma~\ref{lem:maximum-functional}, it is a constrained
critical point as well.
\end{proof}

\begin{proof}[Proof of Theorem~\ref{thm:main} (i)]
Propositions~\ref{prop:local-minimizer} and~\ref{prop:upper-solution}
yield scalar solutions $w^+$ and $w^-$ with negative and positive
energies, respectively. Lemma~\ref{lem:weak-pohozaev} gives their
positivity away from the origin and the multiplier formula
\eqref{eq:multiplier-main}. Proposition~\ref{prop:exact-reduction}
lifts these scalar solutions to synchronized normalized solutions.

The map $w\mapsto(\rho_1w,\rho_2w)$ is a homeomorphism from $\mathcal S$ onto $\mathcal S_{\rm syn}$. By \eqref{eq:energy-factorization}, it transfers the local minimizer on $\mathcal S$ and the minimizer on $\mathcal P^-$ to the corresponding minimizers on $\mathcal S_{\rm syn}$ and $\mathcal P^-_{\rm syn}$, and yields \eqref{eq:opposite-energy-main}.
\end{proof}

\section{The mass-critical and mass-supercritical cases}\label{sec:critical-supercritical}

In this section, we consider the case
\[
q_c\leq q<p.
\]
By Lemma~\ref{lem:critical-supercritical-fibers}, each fiber map admits a unique
critical point, which is the global maximum. Hence the variational analysis is performed on the
whole Pohozaev set $\mathcal P$.

In the mass-critical case, assume that
\eqref{eq:beta-critical-main} holds. Define
\[
m=\inf_{\mathcal P}J .
\]

\begin{lemma}\label{lem:positive-ground-level}
	The variational level satisfies $m>0$. 
\end{lemma}

\begin{proof}
	Suppose first that $q>q_c$. Then $\theta>2$. For any
	$w\in\mathcal P$, the Pohozaev identity,
	Lemma~\ref{lem:weighted-gn} and \eqref{eq:sharp-ckn} imply
	\begin{equation}\label{eq:P-lower-A-super}
		1\leq \delta_q\beta K_q A(w)^{(\theta-2)/2}
		+\kappa S(a,b)^{-p/2}A(w)^{(p-2)/2}.
	\end{equation}
	Since the right-hand side tends to zero as $A(w)\to0^+$,
	there exists a constant $A_0>0$, independent of
	$w\in\mathcal P$, such that $A(w)\geq A_0$.
On $\mathcal P$, using the Pohozaev identity
\[
A(w)=\delta_q\beta B(w)+\kappa C(w),
\]
we obtain
\[
J(w)
=\frac{\theta-2}{2q}\beta B(w)
+\frac dN\kappa C(w)
\geq
c_0A(w)\geq c_0A_0>0,\quad w\in\mathcal P,
\]
where
$
c_0=\min\left\{\frac{\theta-2}{2q\delta_q},\frac dN\right\}>0.
$
Consequently,
$
m=\inf_{\mathcal P}J>0.
$
	
If $q=q_c$, then $\theta=2$. The condition
	\eqref{eq:beta-critical-main} gives
	\[
	c_\beta A(w)\leq \kappa C(w)
	\leq \kappa S(a,b)^{-p/2}A(w)^{p/2}.
	\]
Therefore,
\[
A(w)\geq
\left(\frac{c_\beta S(a,b)^{p/2}}{\kappa}\right)^{\frac{2}{p-2}}
>0,
\qquad w\in\mathcal P .
\]
Moreover, on $\mathcal P$,
\[
J(w)=\frac{d}{N}\kappa C(w)
\geq
\frac{d}{N}c_\beta A(w)>0.
\]
Consequently,
\[
m=\inf_{\mathcal P}J
\geq
\frac dN c_\beta
\left(\frac{c_\beta S(a,b)^{p/2}}{\kappa}\right)^{\frac{2}{p-2}}
>0.
\]
\end{proof}

\begin{lemma}\label{lem:strict-threshold-ground}
We have
\begin{equation*}
 0<m<T_{\kappa}.
\end{equation*}
\end{lemma}

\begin{proof}
    The lower bound follows from Lemma~\ref{lem:positive-ground-level}.
    For the upper bound, let $W$ be the extremal from
    Lemma~\ref{lem:finite-mass-extremal} and set
	\[
	Z=\frac{W}{M(W)^{1/2}}\in\mathcal S .
	\]
	For the pure critical functional, equality in the sharp CKN
	inequality~\eqref{eq:sharp-ckn} gives
	\[
	\max_{t\in\mathbb R}
	\left(
	\frac12e^{2t}A(Z)
	-\frac{\kappa}{p}e^{pt}C(Z)
	\right)
	=
	\frac{d}{N}
	\left(
	\frac{A(Z)}{(\kappa C(Z))^{2/p}}
	\right)^{\frac{N}{2d}}
	=T_\kappa.
	\]
Since
\[
J(t\star Z)
<
\frac12e^{2t}A(Z)-\frac{\kappa}{p}e^{pt}C(Z),
\qquad t\in\mathbb R,
\]
we have
\[
\max_{t\in\mathbb R}J(t\star Z)<T_\kappa .
\]
Let \(t_Z\) be the unique maximizer of the fiber map
\(t\mapsto J(t\star Z)\), whose existence follows from
Lemma~\ref{lem:critical-supercritical-fibers}. Since
\(P(t_Z\star Z)=0\), we have \(t_Z\star Z\in\mathcal P\). Therefore,
\[
m\leq J(t_Z\star Z)
=\max_{t\in\mathbb R}J(t\star Z)
<T_\kappa.
\]
\end{proof}

\begin{proposition}\label{prop:ground-state}
The level $m$ is attained by a nonnegative function $w\in\mathcal P$, which is a constrained critical point of $J|_{\mathcal S}$.  Moreover,
\begin{equation*}
 J(w)=m\in(0,T_{\kappa}).
\end{equation*}
The level admits the minimax characterization
\begin{equation}\label{eq:minimax-characterization}
 m=\inf_{z\in\mathcal S}\max_{t\in\mathbb R}J(t\star z).
\end{equation}
\end{proposition}

\begin{proof}
	We first establish \eqref{eq:minimax-characterization}.
	For every $w\in\mathcal S$, \eqref{eq:dilation-properties}
	and \eqref{eq:fiber-first} give
	$t\star w\in\mathcal S$ and $\Phi_w'(t)=P(t\star w)$
	for all $t\in\mathbb R$.
	By Lemma~\ref{lem:critical-supercritical-fibers},
	$t_*(w)$ is the unique critical point of $\Phi_w$
	and its strict global maximizer.
	Thus $t_*(w)\star w\in\mathcal P$ and
	\[
	m\leq J(t_*(w)\star w)
	=\max_{t\in\mathbb R}J(t\star w).
	\]
	Taking the infimum over $w\in\mathcal S$, we obtain
	\[
	m\leq
	\inf_{w\in\mathcal S}\max_{t\in\mathbb R}J(t\star w).
	\]
	Conversely, for $z\in\mathcal P$,
	$\Phi_z'(0)=P(z)=0$, so the uniqueness of the critical
	point gives $t_*(z)=0$.
	Since $\mathcal P\subset\mathcal S$,
	\[
	\inf_{w\in\mathcal S}\max_{t\in\mathbb R}J(t\star w)
	\leq
	\inf_{z\in\mathcal P}\max_{t\in\mathbb R}J(t\star z)
	=
	\inf_{z\in\mathcal P}J(z)
	=m.
	\]
	Combining the two inequalities proves
	\eqref{eq:minimax-characterization}.
	
	To prove attainment of $m$, let
	$\{z_n\}\subset\mathcal P$ be a minimizing sequence,
	so that $J(z_n)\to m$.
	If $q>q_c$, then $\delta_q/2-1/q>0$, so
	\eqref{eq:energy-on-P-C} and the boundedness of
	$\{J(z_n)\}$ imply that both $\{B(z_n)\}$ and
	$\{C(z_n)\}$ are bounded.
	Since $P(z_n)=0$, we have
	\[
	A(z_n)=\delta_q\beta B(z_n)+\kappa C(z_n),
	\]
	which implies that $\{A(z_n)\}$ is bounded.
	
	If $q=q_c$, then $q\delta_q=2$ and
	$c_\beta=1-\delta_q\beta K_q>0$ by
	\eqref{eq:beta-critical-main}.
	Using \eqref{eq:weighted-gn} and
	\eqref{eq:energy-on-P-C}, together with
	$M(z_n)=1$ and $P(z_n)=0$, we obtain
	\[
	c_\beta A(z_n)
	\leq A(z_n)-\delta_q\beta B(z_n)
	=\kappa C(z_n)
	=\frac{N}{d}J(z_n).
	\]
    Thus, in both cases, $\{z_n\}$ is bounded in $X$, since
    $\|z_n\|_X^2=A(z_n)+1$.
	
Proposition~\ref{prop:full-PS} yields a bounded
Palais--Smale sequence $\{w_n\}\subset\mathcal P$
satisfying \eqref{eq:full-PS-sequence} with $c=m$.
Passing to a subsequence, we may assume that
$w_n\rightharpoonup w$ in $X$.
If critical loss occurs,
Proposition~\ref{prop:compactness-alternative} and
Lemma~\ref{lem:strict-threshold-ground} give
\[
J(w)\leq m-T_\kappa<0.
\]
Thus $w\not\equiv0$. By \eqref{eq:early-mass-recovery},
we have $M(w)=1$, while
Lemma~\ref{lem:weak-pohozaev} gives $P(w)=0$.
Hence $w\in\mathcal P$, and the definition of $m$ yields
\[
0<m\leq J(w)<0,
\]
which is impossible.
Since $m>0$, Proposition~\ref{prop:compactness-alternative}
then implies that, up to a subsequence,
$
w_n\to w \ \text{in }X.
$
By continuity, $w\in\mathcal P$ and
$
J(w)=m.
$
Replacing $w$ by $|w|$ gives a nonnegative minimizer.
By Lemma~\ref{lem:maximum-functional}, it is a constrained
critical point of $J$ on $\mathcal S$.
\end{proof}

\begin{proof}[Proof of Theorem~\ref{thm:main} (ii)]
	Proposition~\ref{prop:ground-state} yields a scalar Pohozaev
	ground state with energy in $(0,T_\kappa)$.
	Lemma~\ref{lem:weak-pohozaev} gives its positivity away from
	the origin and the multiplier formula \eqref{eq:multiplier-main}.
	Proposition~\ref{prop:exact-reduction} lifts this scalar solution
	to a synchronized normalized solution.
	
    The map $w\mapsto(\rho_1w,\rho_2w)$ is a homeomorphism from $\mathcal S$ onto $\mathcal S_{\rm syn}$. By \eqref{eq:energy-factorization}, it transfers the minimizer of $J$ on $\mathcal P$ to a minimizer of $E$ on $\mathcal P_{\rm syn}$, hence to a synchronized Pohozaev ground state. Finally, \eqref{eq:energy-factorization}, together with
    $J(w)=m\in(0,T_\kappa)$, yields
    \eqref{eq:critical-energy-main}.
\end{proof}

\section{Synchronized solutions and the full vector level}\label{sec:main-proof}

\begin{proof}[Proof of Theorem~\ref{thm:weak-coupling}]
	Throughout this proof, $q=q_c$, both prescribed masses are one,
	and $c_\beta=1-\delta_q\beta K_q>0$.
	For $w\in\mathcal S$, set $L(w)=A(w)-\delta_q\beta B(w)$.
	Since $\delta_q=2/q$, we have $L(w)\geq c_\beta A(w)>0$ and
	\begin{equation*}
		J_\chi(t\star w)=\frac12e^{2t}L(w)-\frac\chi p e^{pt}C(w).
	\end{equation*}
	The unique maximum of this fiber is given by
	\begin{equation}\label{eq:critical-coefficient-law}
		\max_tJ_\chi(t\star w)
		=\frac dN\frac{L(w)^{p/(p-2)}}{(\chi C(w))^{2/(p-2)}}.
	\end{equation}
	Taking the infimum over $\mathcal S$ proves
	$m_\chi=\chi^{-2/(p-2)}m_1$ by
	\eqref{eq:minimax-characterization}.
	The positivity and attainment of these levels follow from
	Proposition~\ref{prop:ground-state}.
	The energy factorization for synchronized pairs gives
	$m_{\rm syn}=2m_\kappa$.
	
	For $(u,v)\in\mathcal S_{\rm vec}$, define
	\begin{equation*}
		L_2=L(u)+L(v),\qquad
		K_2=\mu(C(u)+C(v))+p\nu D(u,v).
	\end{equation*}
    Since $L_2,K_2>0$, the fiber map has a unique global
    maximum, given by
	\begin{equation}\label{eq:vector-fiber-maximum}
		\max_tE(t\star u,t\star v)
		=\frac dN\frac{L_2^{p/(p-2)}}{K_2^{2/(p-2)}}.
	\end{equation}
	Each fiber intersects $\mathcal P_{\rm vec}$ at its maximizing
	point. To establish a positive lower bound, set
	$A_2=A(u)+A(v)$. The H\"older inequality gives
	\begin{equation*}
		D(u,v)\leq C(u)^{1/2}C(v)^{1/2}
		\leq\frac12(C(u)+C(v)),
	\end{equation*}
	and hence
	\begin{equation*}
		K_2\leq\kappa(C(u)+C(v))
		\leq\kappa S(a,b)^{-p/2}A_2^{p/2}.
	\end{equation*}
	On $\mathcal P_{\rm vec}$, 
	$K_2=L_2\geq c_\beta A_2$. It follows that
	\begin{equation*}
		A_2\geq\left(\frac{c_\beta S(a,b)^{p/2}}\kappa\right)^{2/(p-2)},
		\qquad E(u,v)=\frac dN K_2\geq\frac dN c_\beta A_2>0.
	\end{equation*}
	Since this lower bound is uniform, $m_{\rm vec}>0$.
	
	Let $U\in\mathcal P_\mu$ be a positive minimizer with $J_\mu(U)=m_\mu$, so that $L(U)=\mu C(U)$. Choose a nonnegative function
	$Z\in C_c^\infty(\mathbb R^N)\cap\mathcal S$ with
	$\operatorname{supp}Z\subset\{x\in\mathbb R^N:1<|x|<2\}$,
	and define $Z_s=s\star Z$ for $s<0$.
	Since $q=q_c$, the scaling identities yield
	\[
		M(Z_s)=1,\quad
		A(Z_s)=e^{2s}A(Z),\quad
		B(Z_s)=e^{2s}B(Z),\quad
		 C(Z_s)=e^{ps}C(Z).
	\]
	Moreover,
	\[
	\operatorname{supp}Z_s
	\subset\{x\in\mathbb R^N:e^{-s}<|x|<2e^{-s}\}.
	\]
	Hence, by the Cauchy--Schwarz inequality and $M(Z_s)=1$,
	\[
	0\leq m_s:=\int_{\mathbb R^N}|x|^{-2a}UZ_s\,dx
	\leq
	\left(
	\int_{|x|>e^{-s}}|x|^{-2a}U^2\,dx
	\right)^{1/2}
	\longrightarrow0
	\qquad\text{as }s\to-\infty,
	\]
	where the convergence follows from
	$U\in L^2(\mathbb R^N;|x|^{-2a})$.
	
	For $0<\varepsilon\leq1/2$, define
	\begin{equation*}
		h_{\varepsilon,s}=\varepsilon U+(1-\varepsilon^2)^{1/2}Z_s,
		\qquad \sigma_{\varepsilon,s}=M(h_{\varepsilon,s})^{1/2},
		\qquad v_{\varepsilon,s}=\frac{h_{\varepsilon,s}}{\sigma_{\varepsilon,s}}.
	\end{equation*}
	Both components of $(U,v_{\varepsilon,s})$ have unit mass.
	Moreover,
	\begin{equation*}
		\sigma_{\varepsilon,s}^2
		=1+2\varepsilon(1-\varepsilon^2)^{1/2}m_s,
		\qquad 1\leq\sigma_{\varepsilon,s},\qquad
		|\sigma_{\varepsilon,s}-1|\leq C m_s.
	\end{equation*}
	In the estimates below, $C$ may depend on the fixed quantities
	$U,Z,\beta,\mu,\nu,p,q$, but is independent of
	$0<\varepsilon\leq1/2$ and $s<0$. Combining the identity
	\begin{equation*}
		v_{\varepsilon,s}-\varepsilon U
		=(\sigma_{\varepsilon,s}^{-1}-1)\varepsilon U
		+\sigma_{\varepsilon,s}^{-1}(1-\varepsilon^2)^{1/2}Z_s
	\end{equation*}
	with the scaling formulas yields
	\begin{align*}
		\||x|^{-a}\nabla(v_{\varepsilon,s}-\varepsilon U)\|_2
		&\leq C(m_s+e^s),\\
		\||x|^{-b}(v_{\varepsilon,s}-\varepsilon U)\|_q
		&\leq C(m_s+e^{2s/q}),\\
		\||x|^{-b}(v_{\varepsilon,s}-\varepsilon U)\|_p
		&\leq C(m_s+e^s).
	\end{align*}
	The preceding estimates yield convergence in the weighted
	gradient and nonlinear norms as $s\to-\infty$.
	Convergence in the weighted $L^2$-norm is precluded by
	$M(v_{\varepsilon,s})=1$ and $M(\varepsilon U)=\varepsilon^2$.
	
	Set $f=|x|^{-b}U$ and
	$g=|x|^{-b}v_{\varepsilon,s}$. Then $f,g\geq0$, and
	$\|g\|_p$ is bounded uniformly for
	$0<\varepsilon\leq1/2$ and $s<0$.
	Using the elementary inequality
	\[
	|\xi^r-\eta^r|
	\leq r(\xi+\eta)^{r-1}|\xi-\eta|,
	\qquad \xi,\eta\geq0,\quad r>1,
	\]
	with $r=p/2$, followed by the H\"older inequality with
	exponents $2$, $2p/(p-2)$, and $p$, we obtain
	\begin{align*}
		\left|D(U,v_{\varepsilon,s})
		-\varepsilon^{p/2}C(U)\right|
		&=
		\left|\int_{\mathbb R^N}
		f^{p/2}\bigl(g^{p/2}-(\varepsilon f)^{p/2}\bigr)\,dx\right|\\
		&\leq
		\frac p2\int_{\mathbb R^N}
		f^{p/2}(g+\varepsilon f)^{p/2-1}
		|g-\varepsilon f|\,dx\\
		&\leq
		\frac p2\|f\|_p^{p/2}
		\bigl(\|g\|_p+\varepsilon\|f\|_p\bigr)^{p/2-1}
		\|g-\varepsilon f\|_p\\
		&\leq
		C\bigl\||x|^{-b}(v_{\varepsilon,s}-\varepsilon U)\bigr\|_p\\
		&\leq C(m_s+e^s).
	\end{align*}
	
	Set
	\begin{align*}
		L_\varepsilon
		&=L(U)+\varepsilon^2A(U)-\delta_q\beta\varepsilon^qB(U),\\
		K_\varepsilon
		&=\mu C(U)(1+\varepsilon^p)+p\nu\varepsilon^{p/2}C(U),\\
		\mathcal E_s
		&=\left(\int_{|x|>e^{-s}}|x|^{-2a}U^2\right)^{1/2}
		+e^s+e^{2s/q}\longrightarrow0.
	\end{align*}
	The preceding norm estimates, together with the same power
	inequality, give
	\begin{equation}\label{eq:vector-coefficient-control}
		|L_2(U,v_{\varepsilon,s})-L_\varepsilon|
		+|K_2(U,v_{\varepsilon,s})-K_\varepsilon|
		\leq C\mathcal E_s.
	\end{equation}
	All four coefficients are uniformly bounded above and bounded
	away from zero, with $L_2,L_\varepsilon\geq L(U)$
	and $K_2,K_\varepsilon\geq\mu C(U)$.
	For $L_\varepsilon$, the lower bound follows from
	$q>2$, $\varepsilon\leq1$, and $L(U)>0$.
	Consequently, the logarithmic and power functions appearing
	in the fiber formulas have bounded derivatives on the
	relevant coefficient ranges. Define
	\begin{equation*}
		t_{\varepsilon,s}
		=\frac1{p-2}\operatorname{log}
		\frac{L_2(U,v_{\varepsilon,s})}{K_2(U,v_{\varepsilon,s})},\quad
		t_\varepsilon=\frac1{p-2}\operatorname{log}\frac{L_\varepsilon}{K_\varepsilon},\quad
		\mathcal F_\varepsilon
		=\frac dN\frac{L_\varepsilon^{p/(p-2)}}{K_\varepsilon^{2/(p-2)}}.
	\end{equation*}
	Thus, the mean value theorem, together with
	\eqref{eq:vector-fiber-maximum} and
	\eqref{eq:vector-coefficient-control}, yields
	\begin{equation}\label{eq:vector-projection-control}
		|t_{\varepsilon,s}-t_\varepsilon|
		+\left|
		\max_{t\in\mathbb R}
		E(t\star U,t\star v_{\varepsilon,s})
		-\mathcal F_\varepsilon
		\right|\leq C\mathcal E_s.
	\end{equation}
	
	Since $0<d<1, N\geq5$, we have $
	2<p=\frac{2N}{N-2d}
	<\frac{2N}{N-2}\leq\frac{10}{3}<4.
	$
	The identity $L(U)=\mu C(U)$ and \eqref{eq:critical-coefficient-law} give
	\[
	\frac{\mathcal F_\varepsilon}{m_\mu}
	=
	\left(
	1+\frac{A(U)}{L(U)}\varepsilon^2
	-\frac{\delta_q\beta B(U)}{L(U)}\varepsilon^q
	\right)^{\frac{p}{p-2}}
	\left(
	1+\frac{p\nu}{\mu}\varepsilon^{p/2}
	+\varepsilon^p
	\right)^{-\frac{2}{p-2}}.
	\]
	Since $q>2$, the first factor equals
	$1+O(\varepsilon^2)$ as $\varepsilon\to0^+$.
	Taylor's formula applied to the second factor gives
	\[
	\left(
	1+\frac{p\nu}{\mu}\varepsilon^{p/2}
	+\varepsilon^p
	\right)^{-\frac{2}{p-2}}
	=
	1-\frac{2p\nu}{(p-2)\mu}\varepsilon^{p/2}
	+O(\varepsilon^p).
	\]
	Combining these expansions and using $p>2$, we find
	constants $C>0$ and $\varepsilon_0\in(0,1/2]$ such that
	\begin{equation}\label{eq:strict-vector-expansion}
		\left|
		\mathcal F_\varepsilon-m_\mu
		+c_3\varepsilon^{p/2}
		\right|
		\leq C\varepsilon^2,
		\qquad
		c_3=\frac{2p\nu}{(p-2)\mu}m_\mu>0,
	\end{equation}
	for every $0<\varepsilon\leq\varepsilon_0$.
	
	Since $p<4$, we may fix
	$\varepsilon\in(0,\varepsilon_0]$ sufficiently small that
	$
	C\varepsilon^{2-p/2}\leq\frac{c_3}{4}.
	$
	With this $\varepsilon$ fixed, the convergence
	$\mathcal E_s\to0$ as $s\to-\infty$ allows us to choose
	$s<0$ sufficiently negative that the error term in
	\eqref{eq:vector-projection-control} satisfies
	$
	C\mathcal E_s\leq\frac{c_3}{4}\varepsilon^{p/2}.
	$
    By construction,
    $
    (t_{\varepsilon,s}\star U,\,
    t_{\varepsilon,s}\star v_{\varepsilon,s})
    \in\mathcal P_{\mathrm{vec}}.
   $
    Combining \eqref{eq:vector-projection-control} with
    \eqref{eq:strict-vector-expansion}, we obtain
    \[
    m_{\mathrm{vec}}
    \leq
    \max_{t\in\mathbb R}
    E(t\star U,t\star v_{\varepsilon,s})
    \leq \mathcal F_\varepsilon+C\mathcal E_s
    \leq m_\mu-\frac{c_3}{2}\varepsilon^{p/2}
    <m_\mu.
    \]
	This establishes the strict comparison
	$m_{\rm vec}<m_\mu$ in
	\eqref{eq:vector-level-comparison} for every $\nu>0$.
    Finally, \eqref{eq:weak-coupling-interval} is equivalent to
    $
    \kappa\leq 2^{(p-2)/2}\mu .
    $
    The scaling relation
    $m_\chi=\chi^{-2/(p-2)}m_1$
    gives
    \[
    m_{\rm vec}<m_\mu
    =\mu^{-2/(p-2)}m_1
    \leq2\kappa^{-2/(p-2)}m_1
    =m_{\rm syn}.
    \]
\end{proof}

\medskip

\noindent{\bf Acknowledgements}

This work has been supported by the National Natural Science Foundation of
China under Grant No. 12371121 and the Fundamental Research Funds for
the Central Universities under Grant No. SWU-KF26002.

\noindent{\bf Data availability statement}

No datasets were generated or analysed during the present study.

\noindent{\bf Conflict of interest}

The authors declare that they have no conflicts of interest.

\end{document}